\documentclass[preprint,12pt]{elsarticle}
\usepackage{amsmath,amssymb,amsthm,txfonts}
 \usepackage{graphicx}
 \usepackage{color}
\usepackage{algorithm,algpseudocode}

\makeatletter
\def\BState{\State\hskip-\ALG@thistlm}
\makeatother
\usepackage{tikz}

\theoremstyle{plain} 
\newtheorem{Theorem}{Theorem} 

\newtheorem{Definition}[Theorem]{Definition}
\newtheorem{definition}[Theorem]{Definition}
\newtheorem{Lemma}[Theorem]{Lemma}
\newtheorem{Proposition}[Theorem]{Proposition}

\theoremstyle{definition}

\theoremstyle{remark}
\newtheorem{examplex}[Theorem]{Example}
\newenvironment{example}
  {\pushQED{\qed}\examplex}
  {\popQED\endexamplex}

\newtheorem{Remark}[Theorem]{Remark}

\def\Hom{\mathop{\rm Hom}\nolimits}

\def\Span{\mathop{\rm Span}\nolimits}

\makeindex

\newcommand{\kT}{\mathcal{T}}
\newcommand{\kN}{\mathcal{N}}

\newcommand{\kA}{\mathcal{A}}
\newcommand{\kM}{\mathcal{M}}
\newcommand{\kF}{\mathcal{F}}

\newcommand{\NN}{\mathbb{N}}

\newcommand{\st}{ such that\;}

\newcommand{\Nf}{\textrm{Nf}}

\newcommand{\cG}{{\sf{G}}}
\newcommand{\cN}{{\sf{N}}}
\newcommand{\cT}{{\sf{T}}}
  
\newcommand{\ck}{{\bf{k}}}

\newcommand{\cB}{{\sf{B}}}

\newcommand{\bM}{{\underline{M}}}

\newcommand{\xdownarrow}[1]{%
  {\left\downarrow\vbox to #1{}\right.\kern-\nulldelimiterspace}
}

\title{Combinatorics of ideals of points.}

\journal{Journal of Algebra}

\begin{document}

\begin{frontmatter}
\title{Combinatorics of ideals of points: \\a Cerlienco-Mureddu-like approach \\ for an iterative lex game.}

\author[label1]{ Michela Ceria}
\author[label2]{Teo Mora}

\address[label1]{Department of Computer Science, University of di Milan michela.ceria@gmail.com}
\address[label2]{DIMA  Universit\`a di Genova theomora@disi.unige.it}

\begin{abstract}
In 1990 Cerlienco and Mureddu gave a combinatorial iterative algorithm which, 
given an ordered set of points, 
returns the lexicographical Gr\"obner escalier of the ideal of these points.
There are many alternatives to this algorithm and in particular, the most efficient
is the Lex Game, which is not iterative on the points, but its performances are
definitely better.
\\
In this paper, we develop an iterative alternative to Lex Game algorithm,
whose performances are very near to those of the original Lex Game, by means of the  Bar Code, 
a diagram which allows to keep track of information on the points and the  corresponding monomials, that are 
lost and usually recomputed many times in Cerlienco-Mureddu algorithm.\\
Using the same Bar Code, we will also give an efficient algorithm to compute squarefree
separator polynomials of the points and the Auzinger-Stetter matrices with respect
to the lexicographical Gr\"obner escalier of the ideal of the points.
\end{abstract}


\end{frontmatter}
\section{Introduction}\label{Introduction}
In 1990 Cerlienco and Mureddu \cite{CeMu, CeMu2, CeMu3} gave a combinatorial algorithm which, 
given an ordered set of points 
$\underline{\mathbf{X}}=[P_1,...,P_N] \subset \ck^n$, $\ck$ a field, 
returns the lexicographical Gr\"obner escalier 
$$\underline{\cN(I({\bf X}))}\subset
\mathcal{T}:=\{x^{\gamma}:=x_1^{\gamma_1}\cdots
x_n^{\gamma_n} \vert \,\gamma:=(\gamma_1,...,\gamma_n)\in \NN^n \}$$ of the vanishing ideal 
$$I({\bf X}):=\{f \in \mathcal{P}:\, f(P_i)=0,\, \forall i\in \{1,...,N\}\}\subset \mathcal{P}:=\mathbf{k}[x_1,...,x_n].$$
Such algorithm actually returns a bijection (labelled 
\emph{Cer\-lien\-co-Mured\-du correspondence} in \cite[II,33.2]{SPES})
$\Phi_{\underline{\mathbf{X}}}:{\underline{\mathbf{X}}} \rightarrow {\underline{\cN(I(\mathbf{X}))}}.$
The algorithm is inductive and thus has complexity $\mathcal{O}\left(n^2N^2\right)$, but it has the advantage of being iterative, in the sense that, given an ordered set of points 
$\underline{\mathbf{X}}=[P_1,...,P_N]$, its related escalier $\underline{\cN(I({\bf X}))}$ 
and correspondence $\Phi_{\underline{\mathbf{X}}}$, for any point $Q\notin\underline{\mathbf{X}}$ it returns a term 
$\tau\in\mathcal{T}$ such that, denoting $\underline{\mathbf{Y}}$ the ordered set $\underline{\mathbf{Y}}:=[P_1,...,P_N,Q]$, 
\begin{itemize}
 \item $\underline{\cN(I({\bf Y}))}=\underline{\cN(I({\bf X}))}\sqcup\{\tau\}$,
 \item $\Phi_{\underline{\mathbf{Y}}}(P_i)=\Phi_{\underline{\mathbf{X}}}(P_i)$ for all $i$ and
 $\tau=\Phi_{\underline{\mathbf{Y}}}(Q)$.
\end{itemize}
\noindent In order to produce the lexicographical Gr\"obner escalier with a better complexity, 
\cite{FRR} gave a completely different approach (\emph{Lex Game}): given a set of 
(not necessarily ordered) points 
$\mathbf{X}=\{P_1,...,P_N\} \subset \ck^n$ they built a trie (\emph{point trie}) 
representing the coordinates of the points and then  used it to build a different trie, the 
\emph{lex trie}, which 
allows to reed the 
lexicographical Gr\"obner escalier $\cN(I({\bf X}))$.
Such algorithm has a very better complexity, $\mathcal{O}\left(nN+N\min(N,nr)\right)$, 
where $r<n$ is the maximal number of edges from a vertex in the point tree, but in order to obtain it, 
\cite{FRR} was forced to give up iterativity.

\smallskip

\noindent In 1982 Buchberger and M\"oller \cite{BM} gave an algorithm (\emph{Buchberger-M\"oller algorithm}) which, for any term-ordering $<$ on $\mathcal{T}$ and any 
set of (not necessarily ordered) points $\mathbf{X}=\{P_1,...,P_N\} \subset \ck^n$ iterating on the $<$-ordered set $\cN(I({\bf X}))$, returns  the Gr\"obner basis of $I({\bf X})$ with respect $<$, 
the set $\cN(I({\bf X}))$ and a family $[f_1,\cdots,f_N]\subset\mathcal{P}$ of separators of $\mathbf{X}$ \emph{id est} a set of polynomials s.t. $f_i(P_j)=\delta_{ij}=\begin{cases}0&i\neq j\cr 
1&i=j.\cr\end{cases}$ 

\noindent Later M\"oller \cite{MMM2} extended the same algorithm to any finite set of functionals defining a 0-dimensional ideal, thus absorbing also the FGLM-algorithm \cite{FGLM} and, on the other 
side, 
proving that Buchberger-M\"oller algorithm has the FGLM-complexity 
\cite{FGLM} $\mathcal{O}(n^2N^3f)$ where $f$ is the avarage cost of evaluating a functional at a term\footnote{A more precise 
evaluation was later given by Lundqvist\cite{Lun}, namely 
$$\mathcal{O}(\min(n,N)N^3+nN^2+nNf+\min(n,N)N^2f).$$}.

\noindent M\"oller \cite{MMM2} gaves also an alternative algorithm (\emph{M\"oller algorithm}) which,
for any term-ordering $<$ on $\mathcal{T}$, given an ordered set of points\footnote{Actually the algorithm 
is stated for an ordered finite set of functionals $[\ell_1,...,\ell_N]\subset\Hom_\ck(\mathcal{P},\ck)$
such that for each $\sigma\leq N$ the set $\{f \in \mathcal{P}:\, \ell_i(f)=0,\, \forall i\leq 
s\}$ is an ideal.} 
$[P_1,...,P_N] \subset \ck^n$,
for each $\sigma\leq N$, denoting $\mathbf{X}_\sigma=\{P_1,...,P_\sigma\}$ returns, with complexity 
$\mathcal{O}(nN^3+fnN^2)$
\begin{itemize}
\item the Gr\"obner basis of the ideal $I({\bf X}_\sigma)$;
\item the correlated escalier $\cN(I({\bf X}_\sigma))$;
\item a term $t_\sigma\in\mathcal{T}$ such that $\cN(I({\bf X}_\sigma)=\cN(I({\bf X}_{\sigma-1}))\sqcup\{\tau\}$,
\item a triangular set $\{q_1,\cdots,q_\sigma\}\subset\mathcal{P}$ s.t. $q_i(P_j)=\begin{cases}0&i< j\cr 1&i=j,\cr\end{cases}$
\item whence a family of separators can be easily deduced by Gaussian reduction,
\item a bijection $\Phi_\sigma$ such that $\Phi_\sigma(P_i)=\tau_i$ for each $i\leq \sigma$, 
which moreover if $<$ is lexicographical, then coincides with Cer\-lien\-co-Mured\-du 
corresondence.
\end{itemize}
Later, Mora \cite[II,29.4]{SPES} remarked that, since the complexity analisis of both Buch\-ber\-ger-M\"ol\-ler and M\"oller algorithm were assuming to perform Gaussian reduction on an $N$-square 
matrix and to evaluate each monomial in the set 
$$\cB(I({\bf X})):=\left\{\tau x_j, \tau\in \cN(I({\bf X}_\sigma)), 1\leq j \leq n\right\}$$ 
over each point $P_i\in\mathbf{X}$, within that 
complexity 
one can use all the information which can be deduced by the computations 
$\tau(P_i),\tau\in\cB(I({\bf X})), 1\leq i\leq N$;
he therefore introduced the notion of \emph{structural description} of a 
0-dimensional ideal \cite[II.29.4.1]{SPES} and gave an algorithm which computes such structural description 
of each ideal $I({\bf X}_\sigma)$.
Also anticipating the recent mood \cite{Mour,Lu,FGLMlike} of degrobnerizing
effective ideal theory, Mora, in connection with Auzinger-Stetter matrices and algorithm \cite{AS} proposed
to present a 
0-dimensional ideal $I\subset \mathcal{P}$ and its quotient algebra $\mathcal{P}/I$ by 
giving its \emph{Gr\"obner representation} \cite[II.29.3.3]{SPES} \emph{id est} the assignement of 
\begin{itemize}
 \item a  $\ck$-linearly independent ordered set $[q_1,\ldots,q_N]\subset\mathcal{P}/I$
 \item $n$ $N$-square matrices $\left(a_{lj}^{(h)}\right), 1\leq h\leq n$,
\end{itemize}
which  satisfiy
\begin{enumerate}
 \item $\mathcal{P}/I\cong\Span_\ck\{q_1,\ldots,q_N\}$,
 \item $x_hq_l=\sum_ja_{lj}^{(h)}q_j, 1\leq j,l\leq N, 1\leq h\leq n.$
\end{enumerate}
\smallskip
\noindent Since M\"oller algorithm and Mora's extension is inductive, 
our aim is to give an algorithm which given an ordered set of points 
$\mathbf{X}=[P_1,...,P_N] \subset \ck^n$ produces for each $\sigma\leq N$
\begin{itemize}
\item the lexicographical Gr\"obner escalier $\cN(I({\bf X}_\sigma))$,
\item the related Cer\-lien\-co-Mured\-du corresondence,
\item a family of squarefree separators for ${\bf X}_\sigma$,
\item the $n$ $N$-square Auzinger-Stetter matrices $\left(a_{lj}^{(h)}\right), 1\leq h\leq n$, which satisfy condition 2. above with respect the linear basis $\cN(I({\bf X}_\sigma))$.
\end{itemize}
The advantage is that, any time a {\em new} point is to be considered, the old data do not 
need to be modified and actually can simplify the computation of the data for the new ideal.
Since the Lex Game approach which has no tool for considering the order of the points has no way of using the data computed for the ideal $I({\bf X}_{\sigma-1})$ in order to deduce those for $I({\bf 
X}_\sigma)$, while M\"oller algorithm and Mora's extension are iterative on the ordered points and intrinsecally produce Cer\-lien\-co-Mured\-du correspondence, in order to achieve our aim, we need to 
obtain a variation of Cer\-lien\-co-Mured\-du algorithm which is not inductive.

\noindent Our tool is the Bar Code \cite{Ce, CeJ}, essentially a reformulation of the point trie which describes in a compact way the combinatorial strucure of a (non necessarily 0-dimensio\-nal) 
ideal; the Bar Code 
allows to remember and reed those data which Cer\-lien\-co-Mured\-du algorithm is forced to inductively recompute. Actually, once the point trie is computed as in \cite{FRR} with inductive complexity 
$\mathcal{O}(N\cdot N\log(N)n).$ the application of the Bar Code allows to compute the lexicographical Gr\"obner escaliers $\cN(I({\bf X}_\sigma))$ and the related 
Cer\-lien\-co-Mured\-du correspondences, with iterative complexity $\mathcal{O}(N\cdot (n+\min(N,nr)))\sim \mathcal{O}(N\cdot nr).$

\noindent 
The families of separators can be iteratively obtain using Lagrange interpolation via data easily deduced from the point trie as suggested in \cite{FRR,Lun} with complexity 
$\mathcal{O}(N\cdot 
\min(N,nr)).$
\\
The computation of the Auzinger-Stetter matrices is based on Lundqvist result \cite[Lemma 3.2]{Lu} and can be inductively performed with complexity\footnote{Naturally, our decision of giving an 
algorithm which can produce data for the the vanishing ideal when a new point is considered 
forbid us of using the new better algorithms for matrix multiplication\cite{BCLRMat,CW,LG,Pan,SV}; thus our complexity 
is  
$\mathcal{O}\left(N^3\right)$ and not $\mathcal{O}\left(N^\omega\right),\omega<2.39$.}
$\mathcal{O}\left(N\cdot (nN^2)\right).$
\\
\smallskip

\noindent After stating the general notation in section \ref{Notations}, we give a brief recap of 
Cerlienco-Mureddu corresponce (section \ref{Cemu}) and the Lex Game algorithm (section \ref{Lex}).
Afterwards (section \ref{BC}), we introduce the Bar Code, that is employed as a tool in our algorithm (section
\ref{core}), whose complexity is discussed in section \ref{complex}.
Finally section \ref{SeparatoriSection} is devoted to separator polynomials and section \ref{AStet} deals with Auzinger-Stetter matrices.
Appendix \ref{PseudoCode} is dedicated to the psudocode of the algorithm, whereas appendix 
\ref{AppA} contains a commented example.

\section{Notations}\label{Notations}
\noindent Throughout this paper we mainly follow the notation of \cite{SPES}. We denote by $\mathcal{P}:=\mathbf{k}[x_1,...,x_n]$ the ring of polynomials in
$n$ variables with coefficients in the field $\ck$. The \emph{semigroup of terms}, generated by the set $\{x_1,...,x_n\}$ is:
$$\mathcal{T}:=\{x^{\gamma}:=x_1^{\gamma_1}\cdots
x_n^{\gamma_n} \vert \,\gamma:=(\gamma_1,...,\gamma_n)\in \NN^n \}.$$
If $t=x_1^{\gamma_1}\cdots x_n^{\gamma_n}$, then $\deg(t)=\sum_{i=1}^n
\gamma_i$ is the \emph{degree} of $t$ and, for each $h\in \{1,...,n\}$
$\deg_h(t):=\gamma_h$ is the $h$-\emph{degree} of $t$.  A \emph{semigroup ordering} $<$ on $\mathcal{T}$  is  a total ordering
\st $ t_1<t_2 \Rightarrow st_1<st_2,\, \forall s,t_1,t_2
\in \mathcal{T}.$ For each semigroup ordering $<$ on $\mathcal{T}$,  we can represent a polynomial
$f\in \mathcal{P}$ as a linear combination of terms arranged w.r.t. $<$, with
coefficients in the base field $\mathbf{k}$:
$$f=\sum_{t \in \mathcal{T}}c(f,t)t=\sum_{i=1}^s c(f,t_i)t_i:\,
c(f,t_i)\in
\mathbf{k}\setminus \{0\},\, t_i\in \mathcal{T},\, t_1>...>t_s,$$ with
$\cT(f):=t_1$   the 
\emph{leading term} of $f$, $Lc(f):=c(f,t_1)$ the  \emph{leading
coefficient} 
of $f$ and $tail(f):=f-c(f,\cT(f))\cT(f)$  the 
\emph{tail} of $f$.
A \emph{term ordering} is a semigroup ordering \st $1$ is lower 
than every variable or, equivalently, it is a \emph{well ordering}.\\
In all paper, we consider the \emph{lexicographical ordering} 
induced
by  $x_1<...<x_n$, i.e:
$$ x_1^{\gamma_1}\cdots x_n^{\gamma_n}<_{Lex} x_1^{\delta_1}\cdots
x_n^{\delta_n} \Leftrightarrow \exists j\, \vert  \,
\gamma_j<\delta_j,\,\gamma_i=\delta_i,\, \forall i>j, $$
which is a term ordering. Since we do not consider any 
term ordering other than Lex, we drop the subscript and denote it by $<$ 
instead of $<_{Lex}$. A subset $J \subseteq \kT$ is a \emph{semigroup ideal} if  $t \in J 
\Rightarrow st \in J,\, \forall s \in \mathcal{T}$; a subset ${\sf N}\subseteq \mathcal{T}$ is an \emph{order ideal} if
$t\in {\sf N} \Rightarrow s \in {\sf N}\, \forall s \vert t$.  
We have that ${\sf N}\subseteq \mathcal{T}$ is an order ideal if and only if 
$\mathcal{T}\setminus {\sf N}=J$ is a semigroup ideal.
\\
Given a semigroup ideal $J\subset\mathcal{T}$  we define ${\sf 
N}(J):=\mathcal{T}\setminus J$. The minimal set of generators ${\sf G}(J)$ of $J$ is called the \emph{monomial basis} 
of $J$. For all subsets $G \subset \mathcal{P}$,  $\cT\{G\}:=\{\cT(g),\, g \in  G\}$ and $\cT(G)$ is the semigroup ideal
of leading terms defined as $\cT(G):=\{t \cT(g),\, t \in \mathcal{T}, g \in G\}$.
Fixed a term order $<$, for any ideal  $I
\triangleleft \mathcal{P}$ the monomial basis of the semigroup ideal 
$\cT(I)=\cT\{I\}$ is called \emph{monomial basis}  of $I$ and denoted again by $\cG(I)$,
whereas the ideal 
$In(I):=(\cT(I))$ is called \emph{initial ideal} and the order ideal 
$\cN(I):=\kT \setminus \cT(I)$ is called \emph{Groebner escalier} of $I$.

\noindent Let $\mathbf{X}=\{P_1,...,P_N\} \subset \ck^n$ be a finite set of distinct points
$$P_i:=(a_{1,i},...,a_{n,i}),\, i=1,...,N.$$
We call
$$I(\mathbf{X}):=\{f \in \mathcal{P}:\, f(P_i)=0,\, \forall i\},$$
the \emph{ideal of points} of $\mathbf{X}$.\\
If we are interested in the \emph{ordered set}, instead of its support $\mathbf{X}$, 
we denote it by  $\underline{\mathbf{X}}=[P_1,...,P_N]$.
The projection maps are defined as follows:\\
\begin{minipage}[b]{0.5\linewidth}
\centering
$$\pi_m:\ck^n \rightarrow \ck^m $$
$$(X_1,..,X_n)\mapsto (X_1,...,X_m), $$
\end{minipage}
\hspace{0.05cm}
\begin{minipage}[b]{0.5\linewidth}
\centering
$$\pi^m:\ck^n \rightarrow \ck^{n-m+1} $$
$$(X_1,..,X_n)\mapsto (X_m,...,X_n)$$
 \end{minipage}
and, for $P \in \ck^n,\; \mathbf{X} \subset \ck^n$, let
$$\Pi_s(P,\mathbf{X}):=\{P_i \in \mathbf{X}  / \pi_s(P_i)=\pi_s(P)\}, $$
$$\Pi^s(P,\mathbf{X}):=\{P_i \in \mathbf{X}  / \pi^s(P_i)=\pi^s(P)\}, $$
extending in the obvious way the meanings of $\pi_s(\alpha),\pi^s(\alpha),\Pi_s(\alpha,D), \Pi^s(\alpha,D)$ to
$\alpha\in  \ck^n$, $D \subset \ck^n.$\\
Moreover, since there is a bijection $\kT \cong \NN^n$ 
between the terms $x_1^{\gamma_1}\cdots x_n^{\gamma_n}$ in $n$ 
variables and the elements $(\gamma_1,...,\gamma_n)$ in $\NN^n$
(i.e. the exponents' lists of the given terms), we 
extend the meaning of the projections $\pi_m,\pi^m$ also to terms, 
writing, by abuse of notation $\pi_m(x_1^{\gamma_1}\cdots x_n^{\gamma_n})
= x_1^{\gamma_1}\cdots x_m^{\gamma_m}$ and $\pi^m(x_1^{\gamma_1}\cdots x_n^{\gamma_n})
=x_1^{\gamma_m}\cdots x_n^{\gamma_n}$.
\\
We recall now some definitions from Graph Theory, following the notation of \cite{FRR}.\\
\begin{Definition}\label{RTree}
We call \emph{tree}  a connected acyclic graph.  A \emph{rooted tree} is a tree where a special \emph{vertex} (or \emph{node})  called \emph{root}  is singled out.
\end{Definition}
\noindent We say that a vertex is on the $h$-th \emph{level} of the tree if its distance from the root is $h$, i.e. 
we have to walk on $h$ edges to come from the root to the given vertex.  If $v$ is a vertex different from the root, 
 and $u$ is the vertex preceding $v$ on the path from the root, then $u$ is the \emph{parent}
 of $v$ and $v$ is a \emph{child} of $u$.
 Two vertices with the same parent are called \emph{siblings}.
 If $v$ is a vertex different from the root and $u$ is on the path
from $v$ to the root, then $u$ is an ancestor of $v$ and $v$ is a descendant of $u$.
 Clearly the root has no parent. We call \emph{leaves} all the vertices
 having no children and we say that a  \emph{branch} is a path from the root to a leaf.
We consider always trees where all branches have the same length.
The vertices lying in the last level of the tree coincide with the  \emph{leaves}; there are no vertices of the tree under them.

\noindent
\begin{Definition}\label{Trie}
 A \emph{trie} is a rooted tree in which there is a symbol written on every edge from a fixed alphabet.
\end{Definition}
Given a trie $T$ we use the following notation:
\begin{itemize}
 \item the $u$-th vertex (read from left to right) at level $h$  is denoted by $v_{h,u}$;
 \item the set of edges connecting $v_{h,u}$ with its children is denoted by $E_{h,u}$;
 \item moreover we associate to each vertex $v_{h,u}$ a set $V_{h,u}$ of indexes\footnote{Since the trie will be associated to a finite set of distinct points {\bf X} (see section \ref{Lex}), the set $V_{h,u}$  has to be thought as a set containing some indexes of points in {\bf X}.} which we call the label of $v_{h,u}$.
\end{itemize}
\begin{example}\label{TrieEsempio}
Consider the trie $T$
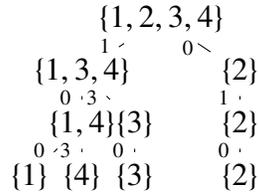
\begin{center}
 \begin{tikzpicture}[>=latex, scale=0.7] 
 \tikzstyle{cerchio}=[circle, draw=black,thick]
 \tikzstyle{cerchior}=[circle, draw=red,thick]
 \tikzstyle{rettangolo}=[rectangle, draw=black,thick]
\node (N_1) at (-0.5,5.75) {$\{1,2,3,4\}$};
\node (N_2) at (-2,4.75) {$\{1,3,4\}$};
\node (N_3) at (-2,3.75) {$\{1,4\}$};
\node (N_4) at (-3,2.75) {$\{1\}$};

\node (N_5) at (1,4.75) {$\{2\}$};
\node (N_6) at (1,3.75) {$\{2\}$};
\node (N_7) at (1,2.75) {$\{2\}$};

\node (N_8) at (-1,3.75) {$\{3\}$};
\node (N_9) at (-1,2.75) {$\{3\}$};
\node (N_10) at (-2,2.75) {$\{4\}$};

\draw [] (N_1) --node[left ]{$\scriptstyle{1}$} (N_2);
\draw [] (N_2) --node[left ]{$\scriptstyle{0}$} (N_3);
\draw [] (N_3) -- node[left ]{$\scriptstyle{0}$}(N_4);

\draw [] (N_1) --node[left ]{$\scriptstyle{0}$} (N_5) ;
\draw [] (N_5) --node[left ]{$\scriptstyle{1}$} (N_6);
\draw [] (N_6) -- node[left ]{$\scriptstyle{0}$}(N_7);

\draw [] (N_2) --node[left ]{$\scriptstyle{3}$} (N_8);
\draw [] (N_8) -- node[left ]{$\scriptstyle{0}$}(N_9);
\draw [] (N_3) -- node[left ]{$\scriptstyle{3}$}(N_10);

\end{tikzpicture}
\end{center}
The vertex  $v_{2,1}$  is labelled by $V_{2,1}=\{1,4\}$; we have $E_{2,1}=\{0,3\}$.
\end{example}

\section{Cerlienco-Mureddu correspondence}\label{Cemu}
In this section, we give a brief description of Cerlienco-Mureddu algorithm,
introduced in \cite{CeMu, CeMu2, CeMu3}, which  is the first combinatorial algorithm that,
given a finite set of distinct points ${\bf X}=\{P_1,...,P_N\}$ computes the lexicographical Groebner
 escalier $\cN(I({\bf X}))$ for the ideal of points of ${\bf X}$.\\
In particular, in \cite{CeMu}, they consider an \emph{ordered} finite set of distinct points in $\ck^n$,
$\underline{{\bf X}}=[P_1,...,P_N]$, and prove that  there  is a one-to-one correspondence between
 $\underline{{\bf X}}$ and the terms of the lexicographical Groebner escalier of $I(\mathbf{X})$:
$$\Phi:\underline{\mathbf{X}} \rightarrow \underline{\cN(I(\mathbf{X}))}$$
$$P_i \mapsto x_1^{\alpha_1^{(i)}}\cdots x_n^{\alpha_n^{(i)}}. $$
They find $\Phi$ using only combinatorics on the coordinates of the elements in ${\bf X}$.
\\
Given $P_i \in {\bf X}$, we denote by $\alpha^{(i)}=(\alpha_1^{(i)},...,\alpha_n^{(i)})\in \ck^n$
the list of exponent of the term $\Phi(P_i)=t_i= x_1^{\alpha_1^{(i)}}\cdots x_n^{\alpha_n^{(i)}}$.
 Sometimes, by abuse of notation, we write\footnote{We recall that there is an isomorphism
between $\kT$ and $\NN^n$.} $\Phi(P_i)=\alpha^{(i)}$.
\\
In order to recall Cerlienco-Mureddu algorithm to compute $\Phi$, we need the following definitions.
\begin{Definition}\label{SigmaValore}
Let ${\bf X}$ be a finite set of distinct points and $P\in \ck^n \setminus {\bf X}$. We call \emph{$\sigma$-value} of $P$ with respect to ${\bf X}$  the maximal integer $s:=\sigma(P,{\bf X}) \in 
\{1,...,n\}$ 
such that $\Pi_{s-1}(P,{\bf X})\neq \emptyset$.
\end{Definition}
\begin{Definition}\label{SigmaAntec}
 Let ${\bf X}$ be a finite set of distinct points, $P\in \ck^n \setminus {\bf X}$
 and $s=\sigma(P,{\bf X})$. We call a point $P_m \in {\bf X}$ \emph{$\sigma$-antecendent} 
 of $P$ w.r.t. ${\bf X}$ and  $\cN(I(\mathbf{X}))$ if $1\leq m\leq \vert{\bf X} \vert=N$ 
 is the maximal integer s.t. 
 \begin{enumerate}
  \item $\pi_{s-1}(P_m)=\pi_{s-1}(P)$
  \item $\pi^{s+1}(\alpha^{(m)})=(0,...,0)$, where $\alpha^{(m)}=\Phi(P_m)$ via Cerlienco-Mureddu correspondence.
 \end{enumerate} 
\end{Definition}
Cerlienco-Mureddu algorithm  is iterative on the points of ${\bf X}$. If $\underline{{\bf X}}=[P_1]$ then $I({\bf X})=(x_1-a_{1,1},...,x_n-a_{1,n})$ and $\cN(I({\bf X}))=\{1\}$.\\
Otherwise, if $N>1$, suppose we already know $\alpha^{(1)}=\Phi(P_1), ..., \alpha^{(N-1)}=\Phi(P_{N-1})$. To find
$\Phi(P_N)=\alpha^{(N)}=(\alpha_1^{(N)},...,\alpha_n^{(N)})$, we have to:
\begin{itemize}
 \item find the $\sigma$-value $s$ of $P_N$ w.r.t. $[P_1,...,P_{N-1}]$;
 \item set $\alpha_n^{(N)}=...=\alpha_{s+1}^{(N)}=0$;
 \item find the $\sigma$-antecedent $P_l$ of $P_N$ w.r.t. $[P_1,...,P_{N-1}]7$
 and $\cN(I(\{P_1,...,P_{N-1}\}))$;
 \item set $\alpha_s^{(N)}=\alpha_s^{(l)}+1$;
 \item find all the points $P_i$ $1 \leq i \leq N$ s.t.$\pi^s(\alpha^{(i)})=(\alpha_s^{(N)}, 0...,0)$ and
 project them w.r.t. the first $s-1$ coordinates finding a set $Y$;
 \item inductively applying Cerlienco-Mureddu algorithm to $Y$,
 compute the Cerlienco-Mureddu correspondence
 $\bar\Phi: Y \rightarrow\cN(I(Y)); $
 \item set $a_i^{(N)} := b_i $ for each $i=1,...,s-1$ where $ (b_1,...,b_{s-1})=\bar\Phi(\pi_{s-1}(P_N))$.
\end{itemize}
There is no complexity analysis in \cite{CeMu}; in \cite{FRR} the authors precise that the number
of comparisons in a straightforward implementation of the algorithm is proportional to $n^2N^2$. 
The quadratic nature is due to the inductive part of the algorithm, which linearly loops over the $n$ variables and the $N$ points.

\section{The Lex Game}\label{Lex}
\noindent The Lex Game algorithm has been introduced by B. Felszeghy, B. R\'{a}th and L. R\'{o}nyai  in \cite{FRR, FR}, as an improvement of
Cerlienco-Mureddu algorithm. In particular, it drops Cerlienco-Mureddu's iterativity, in order to get
the lexicographical Groebner escalier 
with a better complexity, avoiding the squaring produced by the induction. A very precise description of the algorithm, together with a full example and a complexity 
study can be found also in  \cite{Lu}.\\
The first step consists on a preprocessing on the given points, in order to associate them a tree, called \emph{point trie} by Lundqvist.\\
Such a trie, containing the information on the reciprocal relations on the coordinates of the given points is then used to compute a second trie, the \emph{lex trie}, which gives a representation  of the escalier and is actually the solution
to the problem. We see now in details how to construct a point trie, whereas the lex trie computation is only briefly sketched here,
being out of the scope of our paper.
\\
First consider one point $P_1=(a_{1,1},....,a_{n,1})$; its point trie consists of 
only one branch with nodes $v_{0,1},...,v_{n,1}$ s.t. $V_{0,1}=...=V_{n,1}=\{1\}$ and 
$E_{0,1}=\{a_{1,1}\},...,E_{n-1,1}=\{a_{n,1}\}$:
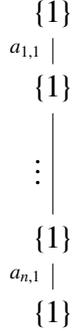
\begin{center}
 \begin{tikzpicture}[>=latex] 
 \tikzstyle{cerchio}=[circle, draw=black,thick]
 \tikzstyle{cerchior}=[circle, draw=red,thick]
 \tikzstyle{rettangolo}=[rectangle, draw=black,thick]
\node (N_1) at (0,5.75) {$\{1\}$};
\node (N_2) at (0,4.75) {$\{1\}$};
\node (N_3) at (0,2.75) {$\{1\}$};
\node (N_4) at (0,1.75) {$\{1\}  $};
\draw [] (N_1) --node[left ]{$\scriptstyle{a_{1,1}}$} (N_2) ;
\draw [] (N_2) --node[left ]{$\scriptstyle{\vdots}$} (N_3);
\draw [] (N_3) -- node[left ]{$\scriptstyle{a_{n,1}}$}(N_4);
\end{tikzpicture}
\end{center}
Suppose now to have constructed the point trie $T:=\mathfrak{T}(\mathbf{X'})$ 
associated to a finite set of distinct
points ${\bf X'}=\{P_1,...,P_{N-1}\}$, and let $P_N=(a_{1,N},....,a_{n,N}) \notin {\bf X'}$.
Adding $P_N$ to the trie means following these steps:
\begin{enumerate}
 \item append to the label of the root $v_{0,1}$ also the index $N$ to the right (i.e.
 modifying the root's label from $V_{0,1}=\{1,...,N-1\}$ to $V_{0,1}=\{1,...,N-1,N\}$); 
 \item for each $1 \leq i \leq n$ consider the children of $v_{i-1,u}$ and the edges in $E_{i-1,u}$, i.e. those
  connecting $v_{i-1,u}$ to its children.
 If one of these edges is labelled by $a_{i,N}$ then let the correspondent child be $v_{i,u'}$; 
 append the index $N$ to the right
 of the label  $V_{i,u'}$ of  $v_{i,u'}$;
 \item when,  at some level $s$, all the edges connected to the children of $v_{s-1,u''}$ are labelled by 
values different from $a_{s,N}$, then consider the rightmost child of  $v_{s-1,u''}$, 
 construct a new node $v_{s,\overline{u}}$ on its right,  labelling it with $V_{s,\overline{u}}=\{N\}$
 and labelling with $a_{i+1,N}$ the edge from $v_{s-1,u''}$ to $v_{s,\overline{u}}$;
 \item for each level $j $ from $s+1$ to $n$ construct a child $v_{j,w}$ of $v_{j-1,z}$ (i.e. the node s.t. $V_{j-1,z}=\{N\}$),  
 labelling it $V_{j,w}=\{N\}$, and label the edge connecting them with $a_{j,N}$.
\end{enumerate}
This way, we get $\mathfrak{T}({\bf X})$, associated to ${\bf X}={\bf X'} \cup \{P_N\}$.
\begin{example}\label{BC Trie}
Given the set $\mathbf{X}=\{(1,0,0),(0,1,0),(1,1,2),(1,0,3)\}$, we display here the construction of 
its point trie $\mathfrak{T}(\mathbf{X})$.
\\
We start with $P_1=(1,0,0)$, associating to it $\mathfrak{T}(\{P_1\})$:
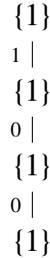
\begin{center}
 \begin{tikzpicture}[>=latex] 
 \tikzstyle{cerchio}=[circle, draw=black,thick]
 \tikzstyle{cerchior}=[circle, draw=red,thick]
 \tikzstyle{rettangolo}=[rectangle, draw=black,thick]
\node (N_1) at (0,5.75) {$\{1\}$};
\node (N_2) at (0,4.75) {$\{1\}$};
\node (N_3) at (0,3.75) {$\{1\}$};
\node (N_4) at (0,2.75) {$\{1\}  $};
\draw [] (N_1) --node[left ]{$\scriptstyle{1}$} (N_2) ;
\draw [] (N_2) --node[left ]{$\scriptstyle{0}$} (N_3);
\draw [] (N_3) -- node[left ]{$\scriptstyle{0}$}(N_4);
\end{tikzpicture}
\end{center}
The second point $P_2=(0,1,0)$ has no common coordinates with $P_1$, so $\mathfrak{T}(\{P_1,P_2\})$ is 
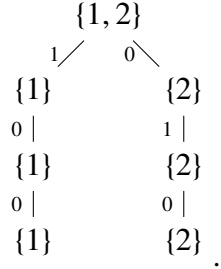
\begin{center}
 \begin{tikzpicture}[>=latex] 
 \tikzstyle{cerchio}=[circle, draw=black,thick]
 \tikzstyle{cerchior}=[circle, draw=red,thick]
 \tikzstyle{rettangolo}=[rectangle, draw=black,thick]
\node (N_1) at (0,5.75) {$\{1,2\}$};
\node (N_2) at (-1,4.75) {$\{1\}$};
\node (N_3) at (-1,3.75) {$\{1\}$};
\node (N_4) at (-1,2.75) {$\{1\}$};

\node (N_5) at (1,4.75) {$\{2\}$};
\node (N_6) at (1,3.75) {$\{2\}$};
\node (N_7) at (1,2.75) {$\{2\}$};

\draw [] (N_1) --node[left ]{$\scriptstyle{1}$} (N_2) ;
\draw [] (N_2) --node[left ]{$\scriptstyle{0}$} (N_3);
\draw [] (N_3) -- node[left ]{$\scriptstyle{0}$}(N_4);

\draw [] (N_1) --node[left ]{$\scriptstyle{0}$} (N_5) ;
\draw [] (N_5) --node[left ]{$\scriptstyle{1}$} (N_6);
\draw [] (N_6) -- node[left ]{$\scriptstyle{0}$}(N_7);
\end{tikzpicture}.
\end{center}
The point $P_3=(1,1,2)$ shares the first coordinate with $P_1$, so for $\mathfrak{T}(\{P_1,P_2,P_3\})$ we get
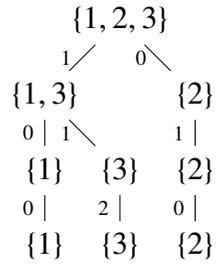
\begin{center}
 \begin{tikzpicture}[>=latex] 
 \tikzstyle{cerchio}=[circle, draw=black,thick]
 \tikzstyle{cerchior}=[circle, draw=red,thick]
 \tikzstyle{rettangolo}=[rectangle, draw=black,thick]
\node (N_1) at (0,5.75) {$\{1,2,3\}$};
\node (N_2) at (-1,4.75) {$\{1,3\}$};
\node (N_3) at (-1,3.75) {$\{1\}$};
\node (N_4) at (-1,2.75) {$\{1\}$};

\node (N_5) at (1,4.75) {$\{2\}$};
\node (N_6) at (1,3.75) {$\{2\}$};
\node (N_7) at (1,2.75) {$\{2\}$};

\node (N_8) at (0,3.75) {$\{3\}$};
\node (N_9) at (0,2.75) {$\{3\}$};

\draw [] (N_1) --node[left ]{$\scriptstyle{1}$} (N_2) ;
\draw [] (N_2) --node[left ]{$\scriptstyle{0}$} (N_3);
\draw [] (N_3) -- node[left ]{$\scriptstyle{0}$}(N_4);

\draw [] (N_1) --node[left ]{$\scriptstyle{0}$} (N_5) ;
\draw [] (N_5) --node[left ]{$\scriptstyle{1}$} (N_6);
\draw [] (N_6) -- node[left ]{$\scriptstyle{0}$}(N_7);

\draw [] (N_2) --node[left ]{$\scriptstyle{1}$} (N_8);
\draw [] (N_8) -- node[left ]{$\scriptstyle{2}$}(N_9);

\end{tikzpicture}.
\end{center}

\noindent The point $P_4=(1,0,3)$ shares the first two coordinates with $P_1$.
The final trie $\mathfrak{T}(\mathbf{X})$ is
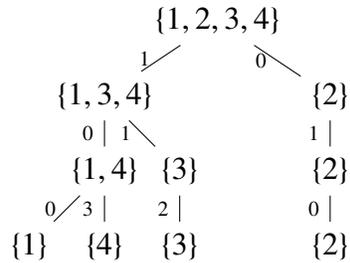
\begin{center}
 \begin{tikzpicture}[>=latex] 
 \tikzstyle{cerchio}=[circle, draw=black,thick]
 \tikzstyle{cerchior}=[circle, draw=red,thick]
 \tikzstyle{rettangolo}=[rectangle, draw=black,thick]
\node (N_1) at (-0.5,5.75) {$\{1,2,3,4\}$};
\node (N_2) at (-2,4.75) {$\{1,3,4\}$};
\node (N_3) at (-2,3.75) {$\{1,4\}$};
\node (N_4) at (-3,2.75) {$\{1\}$};

\node (N_5) at (1,4.75) {$\{2\}$};
\node (N_6) at (1,3.75) {$\{2\}$};
\node (N_7) at (1,2.75) {$\{2\}$};

\node (N_8) at (-1,3.75) {$\{3\}$};
\node (N_9) at (-1,2.75) {$\{3\}$};
\node (N_10) at (-2,2.75) {$\{4\}$};

\draw [] (N_1) --node[left ]{$\scriptstyle{1}$} (N_2);
\draw [] (N_2) --node[left ]{$\scriptstyle{0}$} (N_3);
\draw [] (N_3) -- node[left ]{$\scriptstyle{0}$}(N_4);

\draw [] (N_1) --node[left ]{$\scriptstyle{0}$} (N_5) ;
\draw [] (N_5) --node[left ]{$\scriptstyle{1}$} (N_6);
\draw [] (N_6) -- node[left ]{$\scriptstyle{0}$}(N_7);

\draw [] (N_2) --node[left ]{$\scriptstyle{1}$} (N_8);
\draw [] (N_8) -- node[left ]{$\scriptstyle{2}$}(N_9);
\draw [] (N_3) -- node[left ]{$\scriptstyle{3}$}(N_10);

\end{tikzpicture}
\end{center}
\end{example}
As proved in \cite{Lu}, the computation of the point trie has complexity
 $$O( nN +N \min ( N , nr )),$$ where $r=\max(\vert E_{i,j}\vert)$ denotes the maximal number of edges from a vertex
in the tree and this is also the asymptotic complexity of the whole algorithm (\cite{Lu} Theorem 5.11).
\\
As explained in \cite{FRR}, then, the Lex game takes the point trie and constructs another trie containing
the information on $\cN(I({\bf X}))$. Such a trie is constructed level by level, reading in reversed order the levels of the point trie.
\section{Bar Code for monomial ideals}\label{BC}
\noindent In this section, referring to \cite{Ce}, we summarize the main definitions 
and properties about Bar Codes, which will be used in what follows. First of all, we recall the general definition of  Bar Code. 
\begin{Definition}\label{BCdef1}
A Bar Code $\cB$ is a picture composed by segments, called \emph{bars}, 
superimposed in horizontal rows, which satisfies conditions $a.,b.$ below.
Denote by 
\begin{itemize}
 \item $\cB_j^{(i)}$ the $j$-th bar (from left to right) of the $i$-th row 
 (from top to bottom), i.e. the \emph{$j$-th $i$-bar};
 \item $\mu(i)$ the number of bars of the $i$-th row
 \item $l_1(\cB_j^{(1)}):=1$, $\forall j \in \{1,2,...,\mu(1)\}$ the $(1-)$\emph{length} of the $1$-bars;
 \item $l_i(\cB_j^{(k)})$, $2\leq k \leq n$, $1 \leq i \leq k-1$, $1\leq j \leq \mu(k)$ the $i$-\emph{length} of $\cB_j^{(k)}$, i.e. the number of $i$-bars lying over $\cB_j^{(k)}$
\end{itemize}
\begin{itemize}
 \item[a.] $\forall i,j$, $1 \leq i \leq n-1$, $1\leq j \leq \mu(i)$, $\exists ! \overline{j}\in \{1,...,\mu(i+1)\}$ s.t. $\cB_{\overline{j}}^{(i+1)}$ lies  under  $\cB_j^{(i)}$ 
 \item[b.] $\forall i_1,\,i_2 \in \{1,...,n\}$, $\sum_{j_1=1}^{\mu(i_1)} l_1(\cB_{j_1}^{(i_1)})= \sum_{j_2=1}^{\mu(i_2)} l_1(\cB_{j_2}^{(i_2)})$; we will then say that  \emph{all the rows have the same length}.
\end{itemize}
\end{Definition}
\begin{example}\label{BC1}
 An example of Bar Code $\cB$ is
\\
\begin{minipage}{5cm}
 \begin{center}
\begin{tikzpicture}
\node at (3.8,-0.5) [] {${\scriptscriptstyle 1}$};
\node at (3.8,-1) [] {${\scriptscriptstyle 2}$};
\node at (3.8,-1.5) [] {${\scriptscriptstyle 3}$};

\draw [thick] (4,-0.5) --(4.5,-0.5);
\draw [thick] (5,-0.5) --(5.5,-0.5);
\draw [thick] (6,-0.5) --(6.5,-0.5);
\draw [thick] (7,-0.5) --(7.5,-0.5);
\draw [thick] (8,-0.5) --(8.5,-0.5);
\draw [thick] (4,-1)--(5.5,-1);
\draw [thick] (6,-1) --(6.5,-1);
\draw [thick] (7,-1) --(7.5,-1);
\draw [thick] (8,-1) --(8.5,-1);
\draw [thick] (4,-1.5)--(5.5,-1.5);
\draw [thick] (6,-1.5) --(8.5,-1.5);
\end{tikzpicture}
\end{center}
 
\end{minipage}
\hspace{0.4cm}
\begin{minipage}{8cm}
 $\quad$ \\
The $1$-bars have length $1$.   As regards the other rows, $l_1(\cB_1^{(2)})=2$,
$l_1(\cB_2^{(2)})=l_1(\cB_3^{(2)})= l_1(\cB_4^{(2)})=1$,
$l_2(\cB_1^{(3)})=1$,$l_1(\cB_1^{(3)})=2$ and
 $l_2(\cB_2^{(3)})=$  
\end{minipage}
$l_1(\cB_2^{(3)})=3$, so
 $\sum_{j_1=1}^{\mu(1)} l_1(\cB_{j_1}^{(1)})= \sum_{j_2=1}^{\mu(2)}
l_1(\cB_{j_2}^{(2)})= \sum_{j_3=1}^{\mu(3)} l_1(\cB_{j_3}^{(3)})=5.$
\end{example}
We outline now the construction of the Bar Code associated to a finite set 
of terms.  
First of all, given a  term  $t=x_1^{\gamma_1}\cdots 
x_n^{\gamma_n} \in \mathcal{T}\subset \ck[x_1,...,x_n]$, for each $i \in \{1,...,n\}$, we take
$\pi^i(t):=x_i^{\gamma_i}\cdots x_n^{\gamma_n} \in \mathcal{T}.$ 
Taken a finite set of terms $M\subset \mathcal{T}$, for each $ i \in \{1,...,n\}$, we then define 
$M^{[i]}:=\pi^i(M):=\{s\in \mathcal{T},\,\vert \, \exists t \in M, \pi^i(t)=s\}.$
 \\
Now we take $M\subseteq \mathcal{T}$, with $\vert M\vert =m < \infty$ and we order its 
elements increasingly  w.r.t. Lex, getting the list  
$\bM=[t_1,...,t_m]$. Then, we construct the sets $M^{[i]}$, and 
the corresponding lexicographically ordered lists $\bM^{[i]}$
\footnote{$\bM$ cannot contain repeated terms, while the $\bM^{[i]}$, for $1<i \leq n$, can. 
In case some repeated terms occur in $\bM^{[i]}$, $1<i 
\leq n$, they clearly have to be adjacent in the list, due to the 
lexicographical ordering.}, for $i=1,...,n$.
We can now define the $n\times m $ matrix of terms $\kM$   s.t. 
its $i$-th row is $\bM^{[i]}$, $i=1,...,n$, i.e.
\[\kM:= \left(\begin{array}{cccc}
\pi_{1}(t_1)&... & \pi_{1}(t_m)\\
\pi{2}(t_1)&... & \pi_{2}(t_m)\\
\vdots & \quad &\vdots\\
\pi_{n}(t_1)& ... & \pi_{n}(t_m)
\end{array}\right)\]
\begin{definition}\label{BarCodeDiag}
 The \emph{Bar Code diagram} $\cB$ associated to $M$ (or, equivalently, to 
$\bM$) is a 
$n\times m $ diagram, made by segments s.t. the $i$-th row of $\cB$, $1\leq 
i\leq n$  is constructed as follows:
       \begin{enumerate}
        \item take the $i$-th row of $\kM$, i.e. $\bM^{[i]}$
        \item consider all the sublists of repeated terms, i.e. $[\pi^i(t_{j_1}),\pi^i(t_{j_1 +1}),
        ...,\pi^i(t_{j_1 
+h})]$ s.t. 
        $\pi^i(t_{j_1})= \pi^i(t_{j_1 
+1})=...=\pi^i(t_{j_1 +h})$, noticing that\footnote{Clearly if a term 
$\pi^i(t_{\overline{j}})$ is not 
        repeated in $\bM^{[i]}$, the sublist containing it will be only  
$[\pi_i(t_{\overline{j}})]$, i.e. $h=0$.} $0 \leq h<m$ 
        \item underline each sublist with a segment
        \item delete the terms of $\bM^{[i]}$, leaving only the segments (i.e. 
the \emph{$i$-bars}).
       \end{enumerate}
 We usually label each $1$-bar $\cB_j^{(1)}$, $j \in \{1,...,\mu(1)\}$ with the 
term $t_j \in \bM$.
\end{definition}
\noindent A Bar Code diagram is a 
Bar Code in the sense of definition \ref{BCdef1}.
\begin{example}\label{BarCodeNoOrdId}
Given  $M=\{x_1,x_1^2,x_2x_3,x_1x_2^2x_3,x_2^3x_3\}\subset
\mathbf{k}[x_1,x_2,x_3]$, we have: the $3 \times 5 $ table on the 
left and then to the 
Bar Code on the right:
\\
\begin{minipage}[b]{0.5\linewidth}
\begin{center}
\begin{tikzpicture}
\node at (4.2,-0.5) [] {${\small x_1}$};
\node at (5.2,-0.5) [] {${\small x_1^2}$};
\node at (6.2,-0.5) [] {${\small x_2x_3}$};
\node at (7.2,-0.5) [] {${\small x_1x_2^2x_3}$};
\node at (8.2,-0.5) [] {${\small x_2^3x_3}$};

\node at (4.2,-1) [] {${\small 1}$};
\node at (5.2,-1) [] {${\small 1}$};
\node at (6.2,-1) [] {${\small x_2x_3}$};
\node at (7.2,-1) [] {${\small x_2^2x_3}$};
\node at (8.2,-1) [] {${\small x_2^3x_3}$};

\node at (4.2,-1.5) [] {${\small 1}$};
\node at (5.2,-1.5) [] {${\small1}$};
\node at (6.2,-1.5) [] {${\small x_3}$};
\node at (7.2,-1.5) [] {${\small x_3}$};
\node at (8.2,-1.5) [] {${\small x_3}$};
\end{tikzpicture}
\end{center}
\end{minipage}
\hspace{0.45cm}
\begin{minipage}[b]{0.5\linewidth}
\begin{center}
\begin{tikzpicture}
\node at (4.2,0) [] {${\small x_1}$};
\node at (5.2,0) [] {${\small x_1^2}$};
\node at (6.2,0) [] {${\small x_2x_3}$};
\node at (7.2,0) [] {${\small x_1x_2^2x_3}$};
\node at (8.2,0) [] {${\small x_2^3x_3}$};

\node at (3.8,-0.5) [] {${\scriptscriptstyle 1}$};
\node at (3.8,-1) [] {${\scriptscriptstyle 2}$};
\node at (3.8,-1.5) [] {${\scriptscriptstyle 3}$};

\draw [thick] (4,-0.5) --(4.5,-0.5);
\draw [thick] (5,-0.5) --(5.5,-0.5);
\draw [thick] (6,-0.5) --(6.5,-0.5);
\draw [thick] (7,-0.5) --(7.5,-0.5);
\draw [thick] (8,-0.5) --(8.5,-0.5);
\draw [thick] (4,-1)--(5.5,-1);
\draw [thick] (6,-1) --(6.5,-1);
\draw [thick] (7,-1) --(7.5,-1);
\draw [thick] (8,-1) --(8.5,-1);
\draw [thick] (4,-1.5)--(5.5,-1.5);
\draw [thick] (6,-1.5) --(8.5,-1.5);
\end{tikzpicture}
\end{center}
\end{minipage}
\end{example} 
\noindent Now we recall the vice versa, i.e. how to associate a finite set of terms $M_\cB$ to a given Bar 
Code $\cB$. In \cite{Ce} we first give a more general procedure to do so and then we specialize it in order 
 to have a \emph{unique} set of terms for each Bar Code. Here we give only the specialized version, so 
 we follow the steps below:
\begin{itemize}
 \item[BbC1] consider the $n$-th row, composed by the bars 
$B^{(n)}_1,...,B^{(n)}_{\mu(n)}$. Let $l_1(B^{(n)}_j)=\ell^{(n)}_j$, 
for 
$j\in\{1,...,\mu(n)\}$. Label each bar 
$B^{(n)}_j$ with $\ell^{(n)}_j$ copies 
of $x_n^{j-1}$.
 \item[BbC2] For each $i=1,...,n-1$, $1 \leq j \leq \mu(n-i+1)$ 
 consider the bar $B^{(n-i+1)}_j$ and suppose that it has been 
 labelled by 
$\ell^{(n-i+1)}_j$ copies of a term $t$. Consider all the $(n-i)$-bars 
$B^{(n-i)}_{\overline{j}},...,B^{(n-i)}_{\overline{j}+h}$ 
  lying immediately  above  $ B^{(n-i+1)}_j$; note that $h$ satisfies 
$0\leq h\leq \mu(n-i)-\overline{j}$. 
 Denote the 1-lenghts of 
$B^{(n-i)}_{\overline{j}},...,B^{(n-i)}_{\overline{j}+h}$  by  
$l_1(B^{(n-i)}_{\overline{j}})=\ell^{(n-i)}_{\overline{j}}$,...,
 $l_1(B^{(n-i)}_{\overline{j}+h})=\ell^{(n-i)}_{\overline{j}+h}$. 
 For each $0\leq k\leq h$, label  $ B^{(n-i)}_{\overline{j}+k}$ with 
$\ell^{(n-i)}_{\overline{j}+k}$ copies of $t x_{n-i}^{k}$. 
 \end{itemize}
\begin{Definition}\label{Admiss}
A Bar Code $\cB$ is \emph{admissible} if the set $M$ obtained by applying 
$BbC1$ and $BbC2$ to $\cB$  is an order ideal.
\end{Definition}
\noindent By definition of order 
ideal, using BbC1 and BbC2 is the only way an order 
ideal can be associated to an admissible Bar Code. 
 
\begin{Definition}\label{elist}
 Given a Bar Code $\cB$, 
 let us consider a $1$-bar $B_{j_1}^{(1)}$, with $j_1 
\in \{1,...,\mu(1)\}$.
 The \emph{e-list} associated to $B_{j_1}^{(1)}$ is the $n$-tuple 
$e(B_{j_1}^{(1)}):=(b_{j_1,n},....,b_{j_1,1})$, defined as follows:
 \begin{itemize}
  \item consider the $n$-bar  $B_{j_n}^{(n)}$, lying under 
  $B_{j_1}^{(1)}$. 
The number of $n$-bars on the left of $B_{j_n}^{(n)}$ is  $b_{j_1,n}.$
  \item for each $i=1,...,n-1$, let  $B_{j_{n-i+1}}^{(n-i+1)}$ and 
$B_{j_{n-i}}^{(n-i)}$ be 
  the $(n-i+1)$-bar and the $(n-i)$-bar 
lying under $B_{j_1}^{(1)}$. Consider the $(n-i+1)$-block associated to 
$B_{j_{n-i+1}}^{(n-i+1)}$, i.e. $B_{j_{n-i+1}}^{(n-i+1)}$ and all the bars lying over it. 
The number of $(n-i)$-bars of 
the block, which lie on  the 
left of $B_{j_{n-i}}^{(n-i)}$ is $b_{j_1,n-i}.$
    \end{itemize}
\end{Definition}

\begin{Remark}\label{ElistExp}
 Given a Bar Code $\cB$, fix a $1$-bar  $B_{j}^{(1)}$, with $j \in 
\{1,...,\mu(1)\}$.\\
 Comparing definition \ref{elist} and the steps BbC1 and 
 BbC2 described above, we can observe that the values of the e-list 
$e(B_j^{(1)}):=(b_{j,n},....,b_{j,1})$ are exactly the
exponents of the term 
labelling $B_{j}^{(1)}$, obtained applying BbC1 and BbC2 to
$\cB$.
\end{Remark}

\begin{Proposition}[Admissibility criterion]\label{AdmCrit}
 A Bar Code $\cB$ is admissible if and only if, for each 
 $1$-bar $\cB_{j}^{(1)}$, $j \in \{1,...,\mu(1)\}$, the e-list 
$e(\cB_j^{(1)})=(b_{j,n},....,b_{j,1})$ satisfies the following condition: 
$\forall k \in \{1,...,n\} \textrm{ s.t. } b_{j,k}>0,\, \exists \overline{j} 
 \in \{1,...,\mu(1)\}\setminus \{j\} \textrm{ s.t. } $ $$
e(\cB_{\overline{j}}^{(1)})= (b_{j,n},...,b_{j,k+1}, (b_{j,k})-1, 
b_{j,k-1},...,b_{j,1}). $$
\qed
\end{Proposition}
\noindent Consider the following sets
$$\kA_n:=\{\cB \in \mathcal{B}_n \textrm{ s.t. } \cB  \textrm{ admissible}\} $$
$$\kN_n:=\{\cN \subset \kT,\, \vert \cN\vert < \infty \textrm{ s.t. } \cN 
\textrm{ is an order ideal}\}.$$
We can define the map $\eta: \kA_n \rightarrow \kN_n $; $ \cB \mapsto \cN,$ 
where $\cN$ is the order ideal obtained applying BbC1 and BbC2 to $\cB$,
and it can be easily proved \cite{Ce} that $\eta$ is a bijection.\\
Up to this point, we have discussed the link between Bar Codes and order ideals,
 i.e. we focused on the link between Bar Codes and Groebner escaliers of 
monomial ideals. We show now that, given a Bar Code 
$\cB$ and the order ideal  $\cN =\eta(\cB)$
it is possible to deduce a very specific generating set 
for the monomial ideal $I$ s.t. $\cN(I)=\cN$.
    \begin{Definition}\label{StarSet}
 The \emph{star set} of an order ideal $\cN$ and 
 of its associated Bar Code $\cB=\eta^{-1}(\cN)$ is a set $\kF_\cN$ constructed as 
follows:
 \begin{itemize}
  \item[a)] $\forall 1 \leq i\leq n$, let $t_i$ be a term 
  which labels a $1$-bar lying over $\cB^{(i)}_{\mu(i)}$, 
  then 
  $x_i\pi^i(t_i)\in \kF_\cN$;
  \item[b)] $\forall 1 \leq i\leq n-1$, 
  $\forall 1 \leq j \leq \mu(i)-1$ let 
  $\cB^{(i)}_j$ and $\cB^{(i)}_{j+1}$ be two 
  consecutive bars not lying over the 
same $(i+1)$-bar and let $t^{(i)}_j$ be a term
which labels a $1$-bar lying 
over   $\cB^{(i)}_j$, then 
  $x_i\pi^i(t^{(i)}_j)\in \kF_\cN$.
 \end{itemize}
\end{Definition}
We usually represent $\kF_\cN$ within 
the associated Bar Code $\cB$, inserting
each $t \in \kF_\cN$ on the right of the bar from which 
it is deduced.
Reading the terms from left to right and from the top to
the bottom, $\kF_\cN$ 
is ordered w.r.t. Lex.
\begin{example}\label{BCP}
$\quad $\\
\begin{minipage}{8cm}
For ${\sf
N}=\{1,x_1,x_2,x_3\}\subset
\mathbf{k}[x_1,x_2,x_3]$,  
we have $\kF_\cN=\{x_1^2,x_1x_2,x_2^2,x_1x_3,x_2x_3,x_3^2\}$; looking at 
definition \ref{StarSet}, we can see that  the terms 
\end{minipage}
\hspace{0.5 cm}
\begin{minipage}{5cm}
  \begin{center}
\begin{tikzpicture}[scale=0.4]
\node at (-0.5,4) [] {${\scriptscriptstyle 0}$};
\node at (-0.5,0) [] {${\scriptscriptstyle 3}$};
\node at (-0.5,1.5) [] {${\scriptscriptstyle 2}$};
\node at (-0.5,3) [] {${\scriptscriptstyle 1}$};
 \draw [thick] (0,0) -- (7.9,0);
 \draw [thick] (9,0) -- (10.9,0);
 \node at (11.5,0) [] {${\scriptscriptstyle
x_3^2}$};
 \draw [thick] (0,1.5) -- (4.9,1.5);
 \draw [thick] (6,1.5) -- (7.9,1.5);
 \node at (8.5,1.5) [] {${\scriptscriptstyle
x_2^2}$};
 \draw [thick] (9,1.5) -- (10.9,1.5);
 \node at (11.5,1.5) [] {${\scriptscriptstyle
x_2x_3}$};
 \draw [thick] (0,3.0) -- (1.9,3.0);
 \draw [thick] (3,3.0) -- (4.9,3.0);
 \node at (5.5,3.0) [] {${\scriptscriptstyle
x_1^2}$};
 \draw [thick] (6,3.0) -- (7.9,3.0);
 \node at (8.5,3.0) [] {${\scriptscriptstyle
x_1x_2}$};
 \draw [thick] (9,3.0) -- (10.9,3.0);
 \node at (11.5,3.0) [] {${\scriptscriptstyle
x_1x_3}$};
 \node at (1,4.0) [] {\small $1$};
 \node at (4,4.0) [] {\small $x_1$};
 \node at (7,4.0) [] {\small $x_2$};
 \node at (10,4.0) [] {\small $x_3$};
\end{tikzpicture}
\end{center}
\end{minipage}
$x_1x_3,x_2x_3,x_3^2$ come 
from a), whereas the terms  
$x_1^2,x_1x_2,x_2^2$ come from b).
\end{example}
\noindent In \cite{CMR}, given a monomial ideal $I$, the authors define
 the following set, calling it \emph{star set}:
$$\mathcal{F}(I)=\left\{x^{\gamma} \in \mathcal{T}\setminus {\sf N}(I) \,
\left\vert \,
\frac{x^{\gamma}}{\min(x^{\gamma})} \right. \in {\sf N}(I) \right\}.$$
\begin{Proposition}[\cite{Ce}]\label{DefSt}
With the above notation $\mathcal{F}_{\sf N}=\mathcal{F}(I)$.
\end{Proposition}
\noindent The star set $\kF(I)$ of a monomial ideal $I$ is strongly connected to Janet's 
theory \cite{J1,J2,J3,J4} and to the notion of Pommaret basis \cite{Pom, 
PomAk, SeiB}, as explicitly pointed out in \cite{CMR}.

\section{Our algorithm}\label{core}
\noindent In this section we describe our alternative to Cerlienco-Mureddu algorithm and the Lex game. 
In the next section we will give a complexity analysis and a comparison with
the  aforementioned algoritms. The algorithm that we are going to describe is iterative as Cerlienco-Mureddu algorithm, but it improves the complexity of the execution by noticing that 
Cerlienco-Mureddu algorithm computes \emph{inductively}
some data \emph{more than once}, so the execution time can be reduced and the inductive quadratic complexity can be avoided  by \emph{remembering data in suitable
data structures}. More precisely, we exploit
\begin{itemize}
 \item a point trie $\mathfrak{T}({\bf X})$;
 \item a Bar Code $\cB$;
 \item a matrix $M$.
\end{itemize}
In order to explain the r\^ole of these data structures, we consider a finite set of distinct points
${\bf X}=\{P_1,...,P_N\}\subset \ck^n$ and we order it as $\underline{{\bf X}}=[P_1,...,P_N]$; then
we take $\mathcal{P}=\ck[x_1,...,x_n]$, imposing on it the lexicographical ordering 
induced by $x_1<...<x_n$. The point trie $\mathfrak{T}({\bf X})$, exactly as in the Lex Game, is needed to store the reciprocal relations among the coordinates of the points in 
${\bf X}$. In particular, what we need to know is whether some points share the first $i$ coordinates, for $1\leq i \leq n$. 
As can be deduced by the point trie construction (see Section \ref{Lex}), 
if $1 \leq m,l \leq N$ label the same node at level $1\leq h \leq n$, then
$\pi_h(P_l)=\pi_h(P_m)$.
\begin{example}\label{SameCoord}
 Consider the set ${\bf X}$ and the point trie of example \ref{BC Trie}. We can notice that $P_1=(1,0,0)$ and $P_4=(1,0,1)$ 
 share the first $2$ coordinates, so they label the same nodes of $\mathfrak{T}({\bf X})$ at levels $1,2$, i.e. 
 $1,4 \in V_{1,1}$ and  $1,4 \in V_{2,1}$.
\end{example}
\noindent The Bar Code $\cB$ is the core of our algorithm. The positions of the terms computed in the previous steps are 
exploited to construct, variable by variable, the exponents of the term associated to the point we are dealing with.
The Bar Code is stored in a computer as a list of nested lists; any list corresponds to a bar and its elements are the bars over it; if we are looking at a $1$-bar, its elements are the indexes
of the corresponding points,
since the terms can be desumed by the position (see BbC1-BbC2 in section \ref{BC}).
\begin{example}\label{storeBC}
 Suppose ${\sf
N}=\{1,x_1,x_2,x_3\}\subset
\mathbf{k}[x_1,x_2,x_3]$ is the escalier of a set\\
$ {\bf X}=\{P_1,P_2,P_3,P_4\}$ and $\Phi(P_1)=1,$ 
$\Phi(P_2)=x_2$, $\Phi(P_3)=x_1,$ $\Phi(P_4)=x_3$. The number over the  terms in the $1$-bars identify their corresponding points.
The Bar Code is then stored as $\cB=\Biggl[\biggl[ \Bigl[[1],[3]\Bigr],\Bigl[[2]\Bigr]     \biggr],   \biggl[\Bigl[ [4]    \Bigr]\biggr] \Biggr].$
\end{example}
\noindent Finally we store the terms in the escalier, computed one by one, in a matrix $M$, whose rows are indexed by $[1,...,\vert{\bf X}\vert]$
and whose columns are indexed by $[x_n,...,x_1].$
So for ${\bf X'}=\{P_1,...,P_u\}\subset {\bf X}$ we  have 
\[
M=\begin{bmatrix}
 & \mathbf{x_n} & \mathbf{x_{n-1}} & ... & \mathbf{x_1}\\
 & \mathbf{\downarrow} & \mathbf{\downarrow} & ... & \mathbf{\downarrow}\\
\mathbf{1 \rightarrow } & 0& 0 &... &0\\
\mathbf{2 \rightarrow } & a^{(2)}_n& a^{(2)}_{n-1} &... &a^{(2)}_1\\
\mathbf{...} & ... & ...& ... & ...\\
\mathbf{u \rightarrow } & a^{(u)}_{n}& a^{(u)}_{n-1} &... &a^{(u)}_1
\end{bmatrix}
\]
so that in the entry $(i,j)_{1\leq i \leq u,\, 1 \leq  j \leq n}$ is placed the 
non-negative integer $a^{(i)}_{n-j}$, meaning that the term $t_i \in \cN(I({\bf X'}))$, 
corresponding to $P_i$, has $\deg_{n-j}(t_i)=a^{(i)}_{n-j}$.
We point out that, as can be easily desumed by definition \ref{elist} and Remark \ref{ElistExp}, $M$ contains the 
\emph{e-lists} of the terms in the escalier at the given step, so the values in $M$ also 
allow us to know where the terms are placed in the Bar Code, since the e-list is a sort of set of ``coordinates'' 
of the terms in the Bar Code.
\begin{example}\label{elistExampleExp}
 Consider the Bar Code of example \ref{BCP}. The term $x_2$ has e-list $[0,1,0]$ that, as can be seen in the picture,
 uniquely gives its position in the Bar Code.\\
 \begin{minipage}{8cm}
 Indeed $x_2$ lies over:
 \begin{itemize}
  \item the first $3$-bar $\cB_1^{(3)}$  (indeed, by definition  \ref{elist}, the first $0$ means that there are no $3$-bars on the left  of the $3$-bar of $x_2$.)
 \end{itemize}
 \end{minipage}
\hspace{0.5 cm}
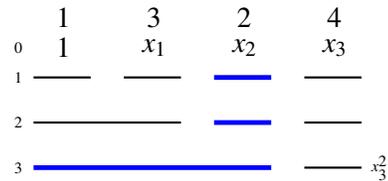
\begin{minipage}{5cm}
  \begin{center}
\begin{tikzpicture}[scale=0.4]
\node at (-0.5,4) [] {${\scriptscriptstyle 0}$};
\node at (-0.5,0) [] {${\scriptscriptstyle 3}$};
\node at (-0.5,1.5) [] {${\scriptscriptstyle 2}$};
\node at (-0.5,3) [] {${\scriptscriptstyle 1}$};
 \draw [line width=1.8pt, color=blue] (0,0) -- (7.9,0);
 \draw [thick] (9,0) -- (10.9,0);
 \node at (11.5,0) [] {${\scriptscriptstyle
x_3^2}$};
 \draw [thick] (0,1.5) -- (4.9,1.5);
 \draw [line width=1.8pt, color=blue] (6,1.5) -- (7.9,1.5);

 \draw [thick] (9,1.5) -- (10.9,1.5);

 \draw [thick] (0,3.0) -- (1.9,3.0);
 \draw [thick] (3,3.0) -- (4.9,3.0);

 \draw [line width=1.8pt, color=blue] (6,3.0) -- (7.9,3.0);

 \draw [thick] (9,3.0) -- (10.9,3.0);

 \node at (1,5.0) [] {\small $1$};
 \node at (1,4.0) [] {\small $1$};
  \node at (4,5.0) [] {\small $3$};
 \node at (4,4.0) [] {\small $x_1$};
 \node at (7,5.0) [] {\small $2$};
 \node at (7,4.0) [] {\small $x_2$};
 
  \node at (10,5.0) [] {\small $4$};

 \node at (10,4.0) [] {\small $x_3$};
\end{tikzpicture}
\end{center}

\end{minipage}
\begin{itemize}
   \item the second $2$-bar over $\cB_1^{(3)}$, i.e. $\cB_2^{(2)}$
  \item the first $1$-bar over $\cB_2^{(2)}$, i.e. $\cB_3^{(1)}$
\end{itemize}

\end{example}
\noindent We start now with the description of our algorithm (Algorithm \ref{IterativeLG}). 
First of all, if we have only one point, i.e. $\vert {\bf X} \vert =N=1$
we set ${\sf N}(1)=\{1\}$ and we construct the point trie $T(P_1)=\mathfrak{T}(\mathbf{X})$ and the 
Bar Code $\cB(1)$ that we display below. The output is stored in the matrix $M$
(lines 2-4).\\
\begin{minipage}{2cm}
\begin{center}
 \begin{tikzpicture}[>=latex, scale=0.7] 
 \tikzstyle{cerchio}=[circle, draw=black,thick]
 \tikzstyle{cerchior}=[circle, draw=red,thick]
 \tikzstyle{rettangolo}=[rectangle, draw=black,thick]
\node (N_1) at (0,5.75) {$\{1\}$};
\node (N_2) at (0,4.75) {$\{1\}$};
\node (N_3) at (0,3.75) {$\{1\}$};
\node (N_4) at (0,2.75) {$\vdots $};
\node (N_5) at (0,1.75)  {$ \{1\} $};
\draw [] (N_1) --node[left ]{$\scriptstyle{a_{11}}$} (N_2) ;
\draw [] (N_2) --node[left ]{$\scriptstyle{a_{21}}$} (N_3);
\draw [] (N_3) -- node[left ]{$\scriptstyle{a_{n-1\,1}}$}(N_4);
\draw [] (N_4) -- node[left ]{$\scriptstyle{a_{n\,1}}$}(N_5);
\end{tikzpicture}
\end{center}
\end{minipage}
\hspace{0.5cm}
\begin{minipage}{4cm}
\begin{center}
\begin{tikzpicture}
\node at (4.7,0.5) [] {${\small 1}$};
\draw [thick] (4.5,0)--(5,0);
\node at (4.7,-0.5)[] {${\small \vdots}$};
\draw [thick] (4.5,-1)--(5,-1);

\node at (3.8,0) [] {${\scriptscriptstyle x_1}$};
\node at (3.8,-1) [] {${\scriptscriptstyle x_n}$};
\end{tikzpicture}
\end{center}
\end{minipage}
\hspace{0.5cm}
\begin{minipage}{3cm}
\[
M=\begin{bmatrix}
 & \mathbf{x_n} & \mathbf{x_{n-1}} & ... & \mathbf{x_1}\\
 & \mathbf{\downarrow} & \mathbf{\downarrow} & ... & \mathbf{\downarrow}\\
\mathbf{1 \rightarrow } & 0& 0 &... &0\\
\end{bmatrix}
\]
\end{minipage}
\\
The above construction  has to be considered as the \emph{base step} for the algorithm.
\\
Now, suppose  $\vert {\bf X} \vert =N>1$ and that the point trie, the matrix and the Bar Code have been constructed for $\{P_1,...,P_{N-1}\}$. We see how to add $P_N$ and get
$\cN(I({\bf X}))$.\\
First, we update the point trie construction, exactly as in the Lex game, by running the subroutine ExTrie (in line 6 Algorithm \ref{ExTrie} is called), which inserts the coordinates of $P_N$. 
In this step we keep track of $s=h$, i.e. the 
level in which $P_N$'s path forks from the pre-existing trie $\mathfrak{T}(\{P_1,...,P_{N-1}\})$. 
By the point trie construction, we can say that the $h$-node $v_{h,u}$ in $P_N$'s path 
(i.e. s.t. $N \in V_{h,u}$) has $V_{h,u}=\{N\}$, whereas the $(h-1)$-node $v_{h-1,u'}$
with $N \in V_{h-1,u'}$ has $\vert V_{h-1,u'} \vert \geq 2$.
This means that such $v_{h-1,u'}$  has at least two children and, by construction, 
the rightmost one, i.e. $v_{h,u}$, is labelled by $V_{h,u}=\{N\}$. \\
Picking the node  $v_{h,u-1}$ on its left and selecting the leftmost element of its label  $V_{h,u-1}$, we get a value $l$: $P_l$ is exactly Cerlienco-Mureddu $\sigma$-antecedent 
(line 14, which calls Algorithm \ref{sant}). Since $P_l$ corresponds to the $l$-th row of $M$, we know the corresponding term  $t_l $ and in particular, this gives us its e-list $e(t_l)$. Using 
$e(t_l)$, we can localize the $s$-bar under $t_l$: 
let it be $\cB_j^{(s)}$. This gives us the following information on $t_N$:
\begin{itemize}
 \item it lies over $\cB_1^{(n)}, \cB_1^{(n-1)},...,\cB_1^{(s+1)}$: indeed, since the $\sigma$-value is $s$, the variables $x_{s+1},...,x_n$ cannot divide the term $t_N$; 
 this implies (by BbC1 - BbC2) that $t_N$ lies over the first 
 $n,...,s+1$ bars, i.e. $a^{(N)}_{s+1}=...=a^{(N)}_n=0$, so $x_n,..., x_{s+1}\nmid t_N$;
 \item it should lie over $\cB_{j+1}^{(s)}$: indeed, similarly,  $\cB_{j}^{(s)}$ is the bar of the $\sigma$-antecedent and 
 $x_s$ must appear in $t_N$ with the exponent incrased by one, w.r.t. that of its  $\sigma$-antecedent, i.e. $a^{(N)}_s=a^{(l)}_{s}+1$.
\end{itemize}
So, we have to test whether  $\cB_{j+1}^{(s)}$ actually lies over $\cB_1^{(n)}, \cB_1^{(n-1)},...,\cB_1^{(s+1)}$; the answer to 
this question (lines 15,16,20) gives rise to two possible cases
\begin{enumerate}
 \item[a.] if $\cB_{j+1}^{(s)}$ does not lie over $\cB_1^{(n)}, \cB_1^{(n-1)},...,\cB_1^{(s+1)}$, then we construct a new $s$-bar of lenght one over  $\cB_1^{(n)}, \cB_1^{(n-1)},...,\cB_1^{(s+1)}$, 
on 
the right of  $\cB_j^{(s)}$,  we clearly label it\footnote{The label depends on the position in the Bar Code, so if we insert a new bar $\cB_{j+1}^{(s)}$, the bars on its right change 
label becoming $\cB_{j+2}^{(s)}$ and so on.} as 
 $\cB_{j+1}^{(s)}$ and we construct a $1,...,s-1$ bar of length $1$ over  $\cB_{j+1}^{(s)}$ (lines 16-19).\\
 Clearly, in this case $t_N=x_s^{j+2}$, so we store the output in the $N$-th row of $M$.
 \item[b.] if $\cB_{j+1}^{(s)}$ lies over $\cB_1^{(n)}, \cB_1^{(n-1)},...,\cB_1^{(s+1)}$, 
we must continue, repeating the procedure, as we   describe below.
\end{enumerate}
In case b., first of all, we have to restrict the point trie only 
to the points whose corresponding terms lie over  $\cB_{j+1}^{(s)}$. 
The set containing these points\footnote{These are the points $P_i$, $1 \leq i < N$ 
s.t. $\pi^{s}(\alpha^{(N)})=\pi^{s}(\alpha^{(i)}) = (\alpha_s^{(N)},0,...,0)$, similarly
to Cerlienco-Mureddu algorithm. The only difference is that we do not consider $P_N$.	} is denoted 
by $S$ and is obtained reading $\cB_{j+1}^{(s)}$. More precisely, $S=\psi(\cB_{j+1}^{(s)})$, where 
$\psi: \cB \rightarrow \kT$ is the function sending each $1$-bar $\cB_l^{(1)}$ in the term $t_l$
over it and, inductively, for $1<u\leq n$,  $\psi(\cB_h^{(u)})=\bigcup_{B\textrm{ over } \cB_h^{(u)} } \psi(B)$ (line 22).
Then, we \emph{read} $P_N$'s path, from level $s-1$ to level $1$, 
looking for the first level $h'$ in which  the node $v_{h',u''}$ with 
 $N \in V_{h',u''}$ is s.t. $S \cap V_{h',u''} \neq \emptyset$.
\\
By the point trie construction, one of the children nodes of $v_{h',u''}$
contains $N$ in its label; we denote it by  $v_{h'+1,w}$.
The new $\sigma$-value is then $s'=h'+1$  and once got it, we can compute the $\sigma$-antecedent exactly as before.  More precisely,
we consider the rightmost child $v_{h'+1,z}$ of $v_{h',u''}$
 on the left of $v_{h'+1,w}$ whose label has nonempty intersection with $S$. The leftmost element in this intersection is the  $\sigma$-antecedent $P_{l'}$.\\
Exploiting
the  e-list of $\Phi(P_{l'})=t_{l'}$, i.e. the $l'$-th row of $M$,  we find the $s'$-bar under, say  $\cB_{j'}^{(s')}$. 
Now, $t_N$ 
\begin{itemize}
 \item lies over $\cB_1^{(n)}, \cB_1^{(n-1)},...,\cB_1^{(s+1)}, \cB_{j+1}^{(s)}, \cB_{e}^{(s-1)},\cB_{e'}^{(s-2)},..., \cB_{e''}^{(s'+1)} $, where 
 $ \cB_{e}^{(s-1)}$ is the first  $(s-1)$-bar over $\cB_{j+1}^{(s)}$, $\cB_{e'}^{(s-2)}$ the first $(s-2)$-bar over $\cB_{e}^{(s-1)}$   and so on;
 \item it should lie over $\cB_{j'+1}^{(s')}$,
\end{itemize}
so we repeat the test on $\cB_{j'+1}^{(s')}$, concluding (case a.) or repeating (case b.) the procedure.
\\
The procedure is repeated until we get to the $1$-bars or if in the 
decision step we get case a.
\\
The position in the Bar Code, being actually the e-list, allows to establish the term 
corresponding to $P_N$, which is finally stored in the $N$-th row of $M$.
\\
For the pseudocode of the described algorithm see Appendix \ref{PseudoCode}.
\begin{Remark}\label{AntecedenteGiusto}
 Let us consider the point trie $\mathfrak{T}(\{P_1,...,P_{N-1}\})$ and suppose we want 
 to deal with $P_N$. 
 We insert $P_N$ in the point trie $\mathfrak{T}(\{P_1,...,P_{N-1}\})$. 
 Let $h$ be the level in which $P_N$'s path
 forks from $\mathfrak{T}(\{P_1,...,P_{N-1}\})$. Clearly the $h$-node 
 $v_{h,u}$ of $P_N$'s path is s.t. $V_{h,u}=\{N\}$ and by the point trie construction, $h$ is the $\sigma$-value.
 \\ Since at level $h$ there is a fork, in the parent node $v_{h-1,u'}$ the label $V_{h-1,u'}$ 
 contains at least one point's index, besides $N$, so 
 $v_{h-1,u'}$ has at least one child besides $v_{h,u}$. 
 Let us take the sibling node
 $v_{h,u-1}$ just on the left of $v_{h,u}$. Suppose that
$v_{h,u-1}$ has label $V_{h,u-1}=\{i_1,...,i_j\}$, $1\leq j\leq N-1$ and consider $i_1$ 
(the leftmost element in the label). We show that $P_{i_1}$ is the point with maximal index
s.t. 
\begin{itemize}
 \item $P_{i_1}$ shares the first $h-1$ coordinates with $P_N$;
 \item the exponents of $x_{h+1},...,x_n$ are zero in $\alpha^{(i_1)}$.
\end{itemize}
First of all, $P_{i_1}$ satisfies these two conditions. Indeed, $i_1,N \in V_{h-1,u}$, so 
$\pi_{h-1}(P_{i_1})=\pi_{h-1}(P_N)$. If 
$P_{i_1}$ would fork from the point trie at some level $k>h$, then there would exist
another point $P_{m_1}$, $m_1<i_1$ with $\pi_{h-1}(P_{m_1})=\pi_{h-1}(P_{i_1})=\pi_{h-1}(P_N)$, 
so also $m_1\in V_{h,u-1}$ being on the left of $i_1$, contradicting the 
choice of the leftmost label. Thus, $P_{i_1}$ forks at some level $\leq h$ and so  
 the exponents of $x_{h+1},...,x_n$ are zero in $\alpha^{(i_1)}$.\\
 Now let us consider some point $P_{n_1}$, with $n_1>i_1$. We show that it cannot satisfy 
the conditions above.  If $n_1\notin V_{h,u-1}$,  it does not share the first $h-1$ coordinates with 
 $P_N$; otherwise it is a label on the right of $i_1$ ($i_1$ was the leftmost). 
 Then,  it must fork from the trie at some level $>h$, this implying that the exponents of 
 $x_{h+1},...,x_n$ are not all zero in $\alpha^{(n_1)}$.
 \end{Remark}
\noindent We show now a commented example of the execution of our algorithm. A more complex one can be found in Appendix \ref{AppA}.
\begin{example}\label{Commented} 
  Consider the following finite set of distinct points
${\bf X}=\{(0,0,0,0), (0,0,0,1), $ $ (0,1,2,3), (1,0,0,0), (1,0,0,1), (1,1,2,3), (0,1,2,4), (1,1,2,4)\}\subset \ck[x_1,x_2,x_3,x_4].$ 
In this case, $n=4$ and we fix on $\ck[x_1,x_2,x_3,x_4]$ the lexicographical ordering induced by $x_1<x_2<x_3<x_4$.
We know that the ideal whose variety is the only point $P_1=(0,0,0,0)$ is $I(\{P_1\})=(x_1,x_2,x_3,x_4)$, so that the 
(lexicographical) Groebner escalier only contains $t_1=1=x_4^0x_3^0x_2^0x_1^0$, so that the matrix 
$M$ of the output monomials' e-lists contains only the zero row:
$M=\begin{bmatrix}
0& 0& 0& 0
\end{bmatrix}
$
We display below both the point trie and the Bar Code.
\\ 
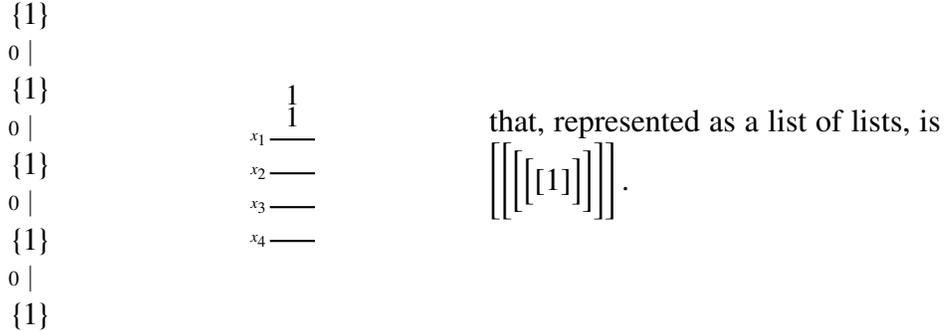
\begin{minipage}{1cm}
 \begin{center}
 \begin{tikzpicture}[>=latex] 
 \tikzstyle{cerchio}=[circle, draw=black,thick]
 \tikzstyle{cerchior}=[circle, draw=red,thick]
 \tikzstyle{rettangolo}=[rectangle, draw=black,thick]
\node (N_1) at (-1,5.75) {$\{1\}$};
\node (N_2) at (-1,4.75) {$\{1\}$};
\node (N_4) at (-1,3.75) {$\{1\} $};
\node (N_9) at (-1,2.75) {$\{1\}$};
\node (N_16) at (-1,1.75) {$\{1\}$};
\draw [] (N_1) --node[left ]{$\scriptstyle{0}$} (N_2) ;
\draw [] (N_2) -- node[left ]{$\scriptstyle{0}$}(N_4);
\draw [] (N_4) --node[left ]{$\scriptstyle{0}$} (N_9);
\draw [] (N_9) -- node[left ]{$\scriptstyle{0}$}(N_16);
\end{tikzpicture}
\end{center}
\end{minipage}
\hspace{0.1cm}
\begin{minipage}{5cm}
\begin{center}
\begin{tikzpicture}[scale=0.3]
\node at (-0.5,-1.5) [] {${\scriptscriptstyle x_4}$};
\node at (-0.5,0) [] {${\scriptscriptstyle x_3}$};
\node at (-0.5,1.5) [] {${\scriptscriptstyle x_2}$};
\node at (-0.5,3) [] {${\scriptscriptstyle x_1}$};
 \draw [thick] (0,-1.5) -- (2,-1.5);
 
 \draw [thick] (0,0) -- (2,0);

 \draw [thick] (0,1.5) -- (2,1.5);

 \draw [thick] (0,3.0) -- (2,3.0);

 \node at (1,4.0) [] {\small $1$};
  \node at (1,5.0) [] {\small $1$};
 \end{tikzpicture}
\end{center}
\end{minipage}
\hspace{0.1cm}
\begin{minipage}{6cm}
that, represented as a list of lists, is $\left[\Biggl[\biggl[ \Bigl[[1] \Bigr]\biggr] \Biggr]\right].$
\end{minipage}
\\
The point $P_2=(0,0,0,1)$, as we can see by the point trie below, forks for $s=4$ and the $\sigma$-antecedent is $P_l,$
for $l=1$, so for now we know that $M[2]=[1,?,?,?]$, i.e. we know only the 
exponent of $x_4$; we have that
$B=\cB_1^{(4)}$ and that there are no next $4$-bars, so that $B'$ still does not exist. We create it getting the  Bar Code on the right:
\\
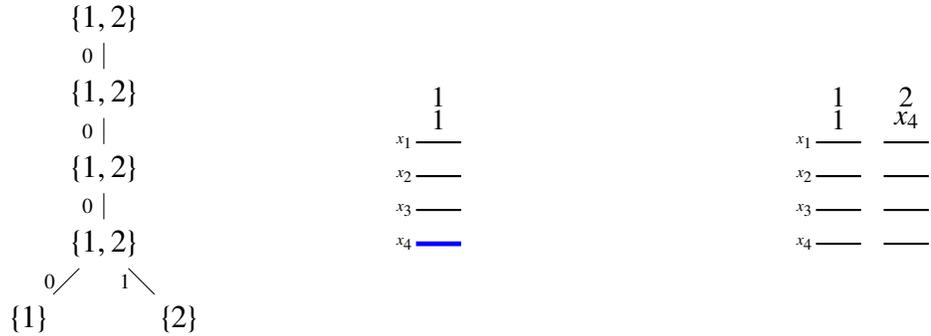
\begin{minipage}{3cm}
\begin{center}

\end{center}
\end{minipage}\\
which, represented by a list, is $\left[ \Biggl[\biggl[ \Bigl[[1] \Bigr]\biggr]\Biggr],\Biggl[\biggl[ \Bigl[[2]\Bigr] \biggr]\Biggr]\right].$
By the Bar Code constructed above, we have  $t_2=x_4$, so
$
M=%
$ 
Take now $P_3=(0,1,2,3)$; after the insertion of this point, the point trie has the following shape:
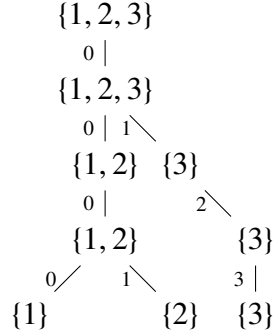
\begin{center}
 %
\end{center}
\noindent The fork happens for $s=2$ and the $\sigma$-antecedent is $P_l,$ for $l=1$ and $B=\cB_1^{(2)}$ so 
$M[3]=[0,0,1,?]$, i.e. $x_3,x_4 \nmid t_3$, $x_2^2 \nmid t_3$ and $x_2 \mid t_3$; we ask for a bar $B'$, lying over $B_1^{(4)},B_1^{(3)}$, 
just on the right of $B$ and since again this bar $B'$ 
does not exist, we create it getting the following Bar Code:\\
\begin{minipage}{4cm}
\begin{center}
%
\end{center}
\end{minipage}
\hspace{0.1cm}
\begin{minipage}{4cm}
which, represented as a list, is 
 $\left[ \Biggl[\biggl[ \Bigl[[1] \Bigr],  \Bigl[[3] \Bigr]  \biggr]\Biggr],\Biggl[\biggl[ \Bigl[[2]\Bigr] \biggr]\Biggr]\right].$
\end{minipage}
\\
Clearly  $t_3=x_2$, so
\[
M=%
\]
For $P_4=(1,0,0,0)$, as shown by the point trie below, the fork happens in $s=1$ and
the $\sigma-$antecedent is $P_l,$ for $l=1$. Then
$M[4]=[0,0,0,1]$. Indeed, we have $B=\cB_1^{(1)}$ 
 and $B'$ has to be created, since there is no $1$-bar on the right of $B$, lying over $\cB_1^{(4)},\cB_1^{(3)},\cB_1^{(2)}$,
so we do it, obtaining the Bar Code on the right
\\
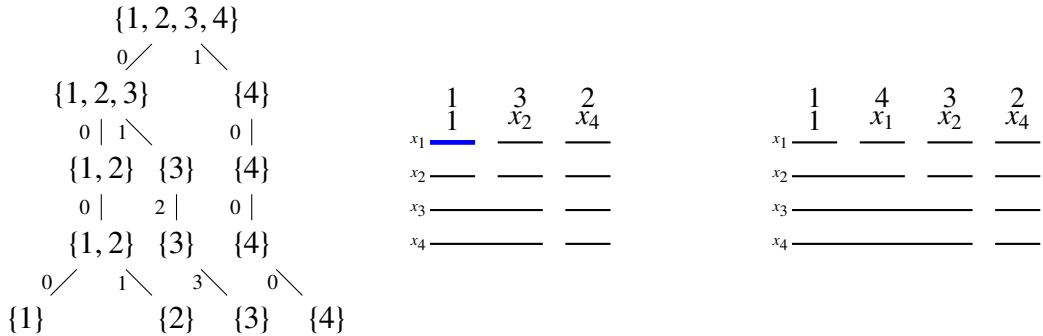
\begin{minipage}{4cm}
\begin{center}
 %
\end{center}
\end{minipage}
\\
 which, represented as a list, is 
 $\left[ \Biggl[\biggl[ \Bigl[[1],[4] \Bigr],  \Bigl[[3] \Bigr]  \biggr]\Biggr],\Biggl[\biggl[ \Bigl[[2]\Bigr] \biggr]\Biggr]\right].$

Clearly  $t_4=x_1$, so
\[
M=%
\]

Consider $P_5=(1,0,0,1)$; the fork happens in $s=4$ and
$l=4$, as shown in the point trie below;

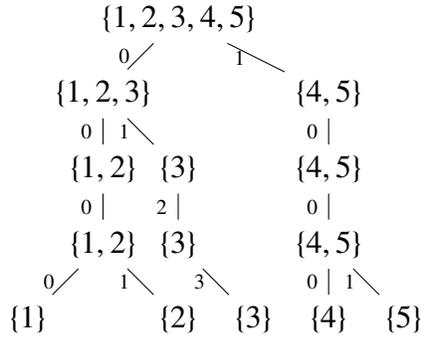
\begin{center}
 %
\end{center}

so $M[5]=[1,?,?,?]$, i.e. $x_4 \mid t_5$.

We have $B=\cB_1^{(4)}$ and in this case $B'=\cB_2^{(4)}$
exists and there is $t_2$ over it, so $S=\{P_2\}$. The fork with $P_2$ happens 
at $s=1$ and the $\sigma$-antecedent is $P_l,$ for $l=2$, so $B=\cB_4^{(1)}$ and $M[5]=[1,0,0,1]$.\\
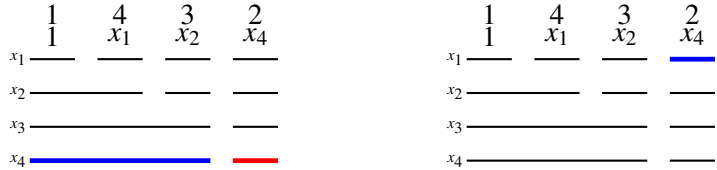
\begin{minipage}{6cm}
\begin{center}
%
\end{center}
\end{minipage}
\\
Since $B'$ still does not exist, we create it and we desume that $t_5=x_1x_4$ and
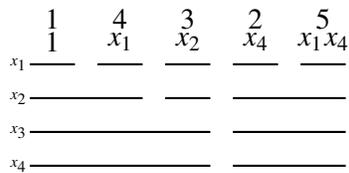
\begin{center}
%
\end{center}
 which, represented as a list, is 
 $\left[ \Biggl[\biggl[ \Bigl[[1],[4] \Bigr],  \Bigl[[3] \Bigr]  \biggr]\Biggr],\Biggl[\biggl[ \Bigl[[2],[5]\Bigr] \biggr]\Biggr]\right].$

So
\[
M=%
\]

After the insertion of  $P_6=(1,1,2,3)$, the point trie has the following shape

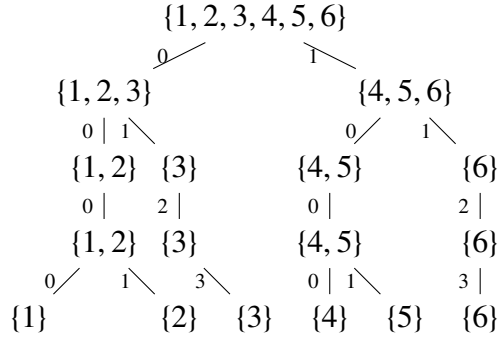
\begin{center}
 %
\end{center}

The fork happens in $s = 2$ and $l = 4$, so  $M[6]=[0,0,1,?]$, $B=\cB_1^{(2)}$ and $B'=\cB_2^{(2)}$, with $t_3=x_2$ over it. 
Then  $S = \{P_3 \}$. The fork with $P_3$ happens in $s=1$ and clearly $l=3$, so $M[6]=[0,0,1,1]$ $B=\cB_3^{(1)}$ and there are no $1$-bars lying over 
$\cB_1^{(4)},\cB_1^{(3)},\cB_1^{(2)}$, so we create $B'$, getting the Bar Code on the right:
\\
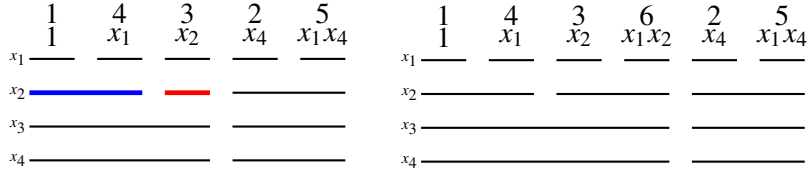
\begin{minipage}{5cm}
\begin{center}
%
\end{center}
\end{minipage}
\\
 which, represented as a list, is 
 $\left[ \Biggl[\biggl[ \Bigl[[1],[4] \Bigr],  \Bigl[[3],[6] \Bigr]  \biggr]\Biggr],\Biggl[\biggl[ \Bigl[[2],[5]\Bigr] \biggr]\Biggr]\right].$

So $t_6=x_1x_2$
\[
M=%
\]

We insert now in the point trie the point $ P_7=(0,1,2,4) $, getting

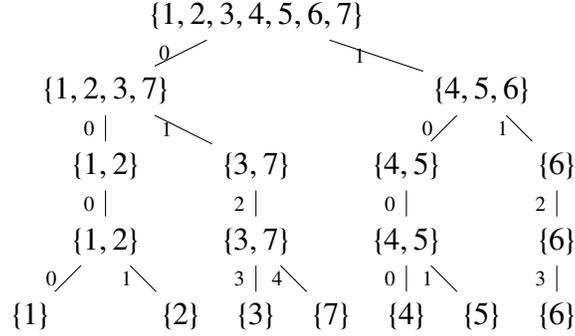
\begin{center}
 %
\end{center}

The fork happens in $s=4$ and the $\sigma$-antecedent is $P_l,$ for $l=3$ so $M[7]=[1,?,?,?]$
We go to $B=\cB_1^{(4)}$ so $B'=\cB_2^{(4)}$, over which we have $t_2,t_5$ so
$S=\{P_2,P_5\}$
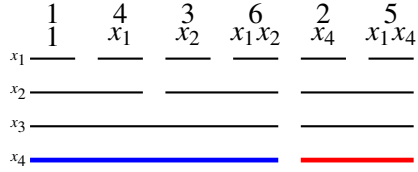
\begin{center}
%
\end{center}

The fork with $S$ happens in $s=2$ and the $\sigma$-antecedent is $P_l,$ for $l=2$; $B=\cB_3^{(2)}$; since there is no $2$-bar 
lying over $\cB_2^{(4)},\cB_2^{(3)}$ and $ B_4^{(2)}$, we have $M[7]=[1,0,1,0]$ and we create $B'$, getting 

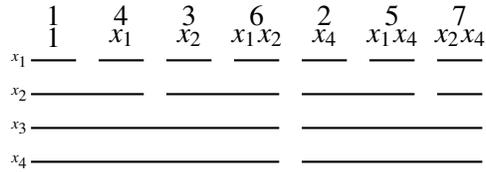
\begin{center}
%
\end{center}
 which, represented as a list, is 
 $\left[ \Biggl[\biggl[ \Bigl[[1],[4] \Bigr],  \Bigl[[3],[6] \Bigr]  \biggr]\Biggr],\Biggl[\biggl[ \Bigl[[2],[5]\Bigr], \Bigl[[7]\Bigr]  \biggr]\Biggr]\right].$\\
The corresponding term is then $t_7=x_2x_4$, so
\[
M=%
\]

Finally we insert $P_8=(1,1,2,4)$ in the point trie and

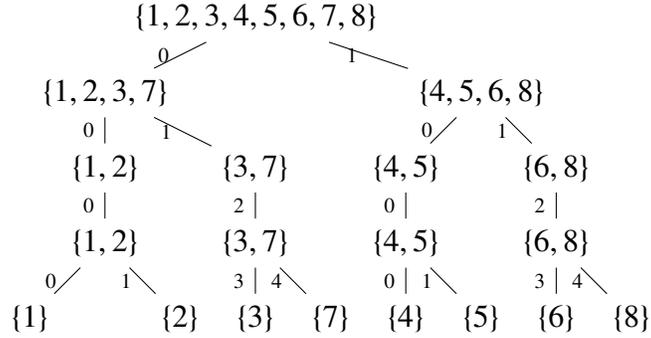
\begin{center}
 %
\end{center}

The fork happens in $s=4$ and the $\sigma$-antecedent is $P_l,$ for $l=6$, so  $M[8]=[1,?,?,?],$ $B=\cB_1^{(4)}$ and  $B'=\cB_2^{(4)}$ with $S=\{P_2,P_5.P_7\}$.

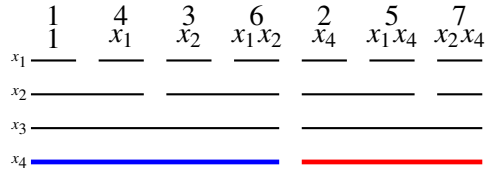
\begin{center}
%
\end{center}

The fork with $S$ happens in $s=2$ and the $\sigma$-antecedent is $P_l,$ for $l=5$ so $M[8]=[1,0,1,?]$,  $B=\cB_3^{(2)}$ and $B'=\cB_4^{(2)}$ and so $S=\{P_7\}$

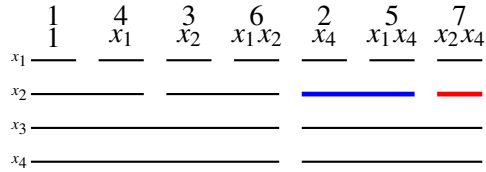
\begin{center}
%
\end{center}

The fork with $S=\{P_7\}$ happens in $s=1$, $l=7$ so $B=\cB_7^{(1)}$. Since there is no $1$ bar over $\cB_2^{(4)},\cB_2^{(3)}, \cB_4^{(2)}$,
then we create it, getting the final Bar Code

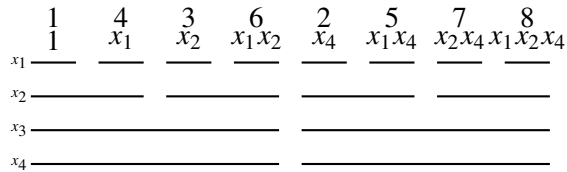
\begin{center}
%
\end{center}
 which, represented as a list, is 
 $\left[ \Biggl[\biggl[ \Bigl[[1],[4] \Bigr],  \Bigl[[3],[6] \Bigr]  \biggr]\Biggr],\Biggl[\biggl[ \Bigl[[2],[5]\Bigr], \Bigl[[7],[8]\Bigr]  \biggr]\Biggr]\right].$

We have $M[8]=[1,0,1,1]$, $t_8=x_1x_2x_4$ and 
finally 

\[
M=%
\]

\noindent We can conclude that the Groebner escalier associated to $I({\bf X})$ is $\cN(I({\bf X}))=\{1,x_1,x_2,x_1x_2,x_4,$
$x_1x_4,x_2x_4, x_1x_2x_4\}.$
\end{example}

\section{Complexity and comparisons}\label{complex}
\noindent In this section, we compute the complexity of our iterative Lex Game algorithm. 
We first remind that  $N$ is the number of points, $n$ is the number of variables and 
$r$ the maximal number of children of a node in the point trie.\\
The construction of the point trie is identical to that of the Lex Game,
with complexity $nN + N \min(N, nr)\sim O(nNr)$.
Interlacing with the Lex Game in order to produce our Bar Code: for each of
the $N$ points, while constructing the new branch of the point trie corresponding to the new point,
we obtain the level $h$ in which the new point forks from the previous ones and the node in the 
$h-1$-level in which the point is still not splitted from the previous ones 
(and in particular the $\sigma$-antecedent). Then 
\begin{itemize}
  \item take the following $h$-block in the Bar Code 
  \item lengthen the $h+1,...,n$ bar under this block and keep track of the corresponding exponents of the monomials (i.e. $0$);
     \begin{enumerate}
       \item if such following  $h$-block in the Bar Code does not exist yet
            we insert it by adding the $1,...,h$ bars (whose length is one)
            
            \item if it exists 
             \begin{itemize}
            \item lengthen also the $h$ bar under this block and keep track of the $h$ exponent of the monomials
             \item we walk in the path of the trie corresponding to the new point, from level $h-1$ to level $1$, repeating the procedure.
            \end{itemize}
     \end{enumerate}

\end{itemize}
The $n$ Bar Code levels (each associated with one of the variables) are read and written each once.
The  cost of detecting the following $h$ block is the same as identifying 
the last element belonging to the current $h$-block and in the $h-1$-node of the trie in which the new point 
appears.
Since the number of points both in the ordered $h$-block and in the $h-1$-ordered node are bounded by $N$, the complexity of this 
problem is $N\log(N)$. Therefore, this procedure costs $N^2n\log(N)$; adding the cost of constructing the point trie, we get again   $N^2n\log(N)$.
We remind that Cerlienco-Mureddu algorithm has complexity $n^2N^2$, while the Lex Game has complexity
 $nN +N \min ( N , nr )$, so our procedure places itself halfway between Cerlienco-Mureddu and the Lex Game, maintaining  Cerlienco-Mureddu's
 iterativity.

 \section{Separator polynomials}\label{SeparatoriSection}
\noindent  In this section, we examine the problem of iteratively computing separator polynomials for 
 a finite set of distinct points.
 \begin{Definition}\label{separators}
 A \emph{family of separators}  for a finite set $\mathbf{X}=\{P_1,...,P_N\}$ of distinct points is a set
 $Q=\{Q_1,....,Q_N\}$ s.t. $Q_i(P_i)=1$ and $Q_i(P_j)=0$, for each $1 \leq i,j \leq N$, $i \neq j$.
 \end{Definition}
 Given  $\mathbf{X}=\{P_1,...,P_N\}$, with $P_i:=(a_{1,i},...,a_{n,i}),\, i=1,...,N,$ we denote by 
 $C=(c_{i,j})$ the witness matrix \cite{Lu}, i.e. the (symmetric) matrix s.t., for $i,j=1,...,N$, $c_{i,j}=0$ if $i=j$ and
 if $i\neq j$, $c_{i,j}= \min\{h: 1 \leq h \leq n \textrm{ s.t. } a_{h,i}\neq a_{h,j}\}$.
 \begin{example}\label{Witnessmat}
 For the set $\mathbf{X}=\{(1,0),(0,1),(0,2)\}$, the witness matrix is 
 \[
C=\begin{bmatrix}
0& 1& 1\\
1& 0 & 2\\
1& 2 & 0
\end{bmatrix}
\]
 \end{example}
We compute the separators iteratively on the points, by means of a variation of the following Lagrange formula:
$$R_i=\prod_{i\neq j} \frac{x_{c_{i,j}}- a_{c_{i,j},j}}{a_{c_{i,j},i}-a_{c_{i,j},j}} =
\prod_{j \neq i} p_{i,j}^{[c_{i,j}]},  \textrm{ with }p_{i,j}^{[c_{i,j}]} =  
\frac{x_{c_{i,j}}- a_{c_{i,j},j}}{a_{c_{i,j},i}-a_{c_{i,j},j}}.$$
In particular, by means of the point trie, we can get squarefree separator polynomials $Q_i=\sqrt{R_i}$, as in \cite{FRR,MMM2}, 
differently from \cite{Lun}, even if in \cite[Remark 5.3]{Lun}, the author remarks that removing exponents to all factors of separator polynomials, produces again separator polynomials. We have to point out also that our polynomials have also smaller degree than those in 
\cite{FRR}, even if also the polynomials found there are squarefree.\\
If there is only one point in $\mathbf{X}$, then we set $Q_1=1$.\\
Suppose now to have computed the separators $Q_1,...,Q_{N-1}$ associated to $\{P_1,...,P_{N-1}\}$ and 
to deal with $P_N$. We see now how to compute the separator of $P_N$ and how to modify the previous ones, to get 
the new separators  $Q_1',...,Q_{N}'$ for $\mathbf{X}$.
\\
First we set $Q_N'=1$.
Then, for each variable $j=1,...,n$, we take the node $v_{j,u}$ s.t. $N \in V_{j,u}$.
For each sibling $v_{j,u'}$ of  $v_{j,u}$, we pick an element $\overline{i} \in  V_{j,u'}$
and we set $Q_N'=Q_N'p_{N,\overline{i}}^{[j]}$.
If $\vert V_{j,u}\vert =1$, namely $v_{j,u}$ is labelled only by $N$, then, for each 
sibling $v_{j,u'}$, for each $i\in  V_{j,u'}$, we set $Q_i'=Q_ip_{i,N}^{[j]}$.\\
Once concluded this procedure, if a separator $Q_h$, $\leq h \leq N$ has \emph{not} been involved
in the above steps, we set $Q_h'=Q_h$, getting a family of separators $\{Q_1',...,Q_N'\}$ 
for $\mathbf{X}=\{P_1,...,P_N\}$.\\
The complexity of a single iterative round of our procedure is $\mathcal{O}(\min(N,nr)).$
\\In the following simple example, we compare the separator families one gets using our method and those 
stated in \cite{FRR,MMM2,Lun}
\begin{example}\label{separatoriconfrontati}
 Consider the set $\mathbf{X}=\{P_1=(1,0),P_2=(0,1),P_3=(0,2)\}$ of example \ref{Witnessmat}. 
 The formula stated in \cite{FRR}, gives the 
 following separator family: 
 $$Q_1=\frac{1}{2}x(y-1)(y-2); \, Q_2=y(x-1)(y-2);\, Q_3=-\frac{1}{2}(x-1)y(y-1),$$
 whereas, using the formula in \cite{Lun}, we get
 $$Q_1=x^2;\, Q_2=(x-1)(y-2);\, Q_3=-(x-1)(y-1).$$ 
 The separator polynomials we get from Moeller algorithm \cite{MMM2} by trivial reduction on the triangular
 polynomials are $Q_1=x$, $Q_2=2-2x-y$ and $Q_3=x+y-1$; 
 these polynomials are squarefree and with support in the escalier associated to $I(\mathbf{X})$.
\\
 The point trie we can construct from our set is 
  \begin{center}
 \begin{tikzpicture}[>=latex, scale=0.7] 
 \tikzstyle{cerchio}=[circle, draw=black,thick]
 \tikzstyle{cerchior}=[circle, draw=red,thick]
 \tikzstyle{rettangolo}=[rectangle, draw=black,thick]
\node (N_1) at (-0.5,5.75) {$\{1,2,3\}$};
\node (N_2) at (-2,4.75) {$\{1\}$};
\node (N_3) at (-2,3.75) {$\{1\}$};

\node (N_5) at (1,4.75) {$\{2,3\}$};
\node (N_6) at (0,3.75) {$\{2\}$};
 \node (N_7) at (2,3.75) {$\{3\}$};
\draw [] (N_1) --node[left ]{$\scriptstyle{1}$} (N_2);
\draw [] (N_2) --node[left ]{$\scriptstyle{0}$} (N_3);
\draw [] (N_1) --node[left ]{$\scriptstyle{0}$} (N_5) ;
\draw [] (N_5) --node[left ]{$\scriptstyle{1}$} (N_6);
\draw [] (N_5) -- node[left ]{$\scriptstyle{2}$}(N_7);
\end{tikzpicture}
\end{center}
In the first step, we set $Q_1''=1$; then, adding $P_2$ to the trie we set $Q_2'=p^{[1]}_{2,1}=-(x-1)$ and we modify 
also $Q_1''$, setting $Q_1'=Q_1''p^{[1]}_{1,2}=x$, since, when $P_3$ is still not in the trie, the node $v_{1,2}$,
has $V_{1,2}=\{2\}$. So, w.r.t. $\{P_1,P_2\}$, we have $Q_1'=x$, $Q_2'=-(x-1)$. Finally, we add $P_3$. This
way, $Q_3=p^{[1]}_{3,1}p^{[2]}_{3,2}=-(x-1)(y-1)$ and since $V_{2,3}=\{3\}$, $Q_2=Q_2'p^{[2]}_{2,3}=(x-1)(y-2)$.
Finally, we have 
$$Q_1=x;\, Q_2=(x-1)(y-2);\, Q_3=-(x-1)(y-1), $$
which are exactly the polynomials that one could obtain from those of 
\cite{Lun} by removing exponents from each factor.
\end{example}
 
\section{Auzinger-Stetter matrices}\label{AStet}
\noindent In this section, we finally deal with the computation of Auzinger-Stetter
matrices associated 
to the zerodimensional (radical) ideal of a finite set of distinct points $\mathbf{X}=\{P_1,...,P_N\}$.
\begin{Definition}\label{ASMat}
 Let $I\triangleleft \ck[x_1,...,x_n]$ be a zerodimensional ideal and
 $A:=\ck[x_1,...,x_n]/I$. For each $f \in A$ we can denote  $\Phi_f:A\rightarrow A$ the linear
 form describing the multiplication by $f$ in $A$
 and, fixed a basis\footnote{where we denote $[f]\in A$ the residue class modulo $I$ 
 of an element $f\in \ck[x_1,...,x_n]$.} $B=\{[b_1],\ldots,[b_m]\}$ for $A$, we can represent it by a matrix $A_f=(a_{ij})$
 so that
$[b_if]=\sum_j   a_{ij} [b_j]$ for each $i$.
 \\ We call \emph{Auzinger-Stetter matrices} associated to $I$, the 
 matrices $A_{x_i},$ for $i=1,...,n$, defined with respect to the basis given by the lex escalier 
of I.
\end{Definition}
\noindent Given  finite set of distinct points
$\mathbf{X}=\{P_1,...,P_N\}$, we denote by $I:=I(\mathbf{X})\triangleleft \ck[x_1,...,x_n]$ 
the associated zerodimensional radical ideal. If $H=\{f_1,...,f_m\}\subset \ck[x_1,...,x_n]$ is a set of polynomials,
we denote by $H(\mathbf{X})=(h_{i,j})$ the $m\times N$ matrix s.t. $h_{i,j}=f_i(P_j)$,
$1 \leq i \leq m$, $1 \leq j \leq N$. Let us denote as usual by $\cN(I):=\{t_1,...,t_N\}$ the lexicographical Groebner 
escalier associated to $I$; given a term $s \in \kT$, $s\cN(I):=\{st_1,...,st_N\}$ and    
\begin{itemize}
\item $A_{x_h}:=\left(a_{li}^{(h)}\right)_{li}, 1\leq h\leq n$, $1\leq l,i\leq N$, the  Auzinger-Stetter matrices w.r.t. $\cN(I)$;
\item $B:=\cN(I)(\mathbf{X}):= \left(b_{lj}\right)_{lj}$, $1\leq l,j\leq N$, $b_{lj}:=t_l(P_j)$;
\item $C:= \left(c_{ji}\right)_{ji}$, $1\leq j,i\leq N$, the inverse matrix of $B$, i.e.
$C := B^{-1}$
\item $D^{(h)}:=\left(d_{lj}^{(h)}\right)_{lj}, 1\leq h\leq n$, $1\leq l,j\leq N$, $d^{(h)}_{lj}:=\alpha_h^{(j)}t_l(P_j)$,
the evaluation of $x_ht_l$ at the point $P_j$.
\end{itemize}
Our tool is the result proposed in the following
\begin{Lemma}[\cite{Lu}, Lemma 3.2]
 Let $\mathbf{X}=\{P_1,...,P_N\}$ be a finite set of distinct points and $I:=I(\mathbf{X})\triangleleft \ck[x_1,...,x_n]$ 
 the associated zerodimensional radical ideal.
 \\
 Let $\cN=\{t_1,...,t_N\}\subset \ck[x_1,...,x_n]$ such that $[\cN]=\{[t_1],...,[t_N]\}$ is a basis for $A:=\ck[x_1,...,x_n]/I$. 
 Then, for each $f \in \ck[x_1,...,x_n]$ we 
have
 $$\Nf(f,\cN)=(t_1,...,t_N)(\cN(\mathbf{X})^{-1})^t (f(P_1),...,f(P_N))^t,$$
 where $\Nf(f,\cN)$ is the normal form of $f$ w.r.t. $\cN$, i.e. the unique expression of the residue modulo $I$ of $f$ as 
 a linear combination of the elements in the basis $[\cN]$.
\end{Lemma}
In particular, to compute $A_{x_h}$, we can employ the above lemma as follows. \\
We first point out that, for $1 \leq l \leq N$, the $l$-th 
row of $A_{x_h}$ is the normal form of $x_ht_l$:
$$\Nf(x_ht_l,\cN(I))=\sum_{i=1}^N a_{li}t_i=(t_1,...,t_N)C^t(x_ht_l(P_1),...,x_ht_l(P_N))^t=$$ 
$$(t_1,...,t_N)C^t
(d_{l1}^{(h)},...,d_{lN}^{(h)})^t=\sum_i \left(\sum_{j=1}^N d_{lj}^{(h)}c_{ji}\right) t_i.$$
This trivially implies that 
$$A_{x_h}= D^{(h)}C.$$
Of course we need a procedure which iteratively extends the computation of the inverse of an invertible matrix;
we cannot use the recent fast invertible matrix approach \cite{BCLRMat,CW,LG,Pan,SV}, thus we only aim to obtain a $\mathcal{O}\left(N\cdot (nN^2)\right)$ complexity and we limit ourselves to an 
adaptation of the computation of the inverse of a $N\times N$ matrix
by Gaussian reduction on the columns.\\
In particular  we begin writing below the $N\times N$ matrix $B$ an identical $N$-square matrix $I$, getting $\begin{pmatrix}B\\ \cline{1-1} I\\ \end{pmatrix}$. Then, we perform
column reduction to the whole $\begin{pmatrix}B\\ \cline{1-1} I\\ \end{pmatrix}$, in order to transform $B$ to the identity. 
In particular, at each step, one gets a new matrix $\begin{pmatrix}F\\ \cline{1-1} E\\ \end{pmatrix}$,
where clearly $F=(f_{li}),\, E=(e_{l,i})$ are   $N$-square matrices, s.t. the following relations on
columns hold
\begin{equation}\label{fBe}
 (f_{1i}...f_{Ni})^t=B(e_{1i}...e_{Ni})^t,\, 1\leq i \leq N.
\end{equation}
Each reduction step, clearly modifies also $E$, so once concluded reduction on $B$
we obtain a matrix $\begin{pmatrix}I\\ \cline{1-1} C\\ \end{pmatrix}$ and we have $C=B^{-1}$.\\
Now we adapt the strategy above in order to make the Auzinger-Stetter matrices computation iterative, 
so we will compute the matrices $B,C,D^{(h)}$
and finally $A_{x_h}$, $1\leq h \leq n$, iteratively on the points of $ \mathbf{X}$.
\\
If $\mathbf{X}=\{P_1\}$, then $\cN(I(\mathbf{X}))=\{1\}$, $B=C=1$ and, if $P_1=(a_{1,1},...,a_{n,1})$, then $D^{(h)}=(a_{h,1})=A_{x_h}$.\\
Now, if $\vert \mathbf{X} \vert >1$, we can construct 
the Auzinger-Stetter matrices iteratively, exploiting the Bar Code and again 
our iterative alternative to Lex Game algorithm.
Suppose we have $\mathbf{X}'=\{P_1,...,P_{N-1}\}$ and we have computed 
$\cN(I(\mathbf{X}'))$ via our algorithm, 
together with the associated Bar Code $\cB'$; we add $P_N$, getting $\mathbf{X}=\{P_1,...,P_{N}\}$.
We first point out that, in our Lex Game alternative algorithm, when we conclude 
the iterative step for $P_N$, we take the last $\sigma$-value $s$ and 
the corresponding $\sigma$-antecedent $P_l$ and we set $t_N=x_st_l$, then, in the 
construction of the matrix $B$, we can compute iteratively the elements $b_{N,i}=t_N(P_i)$, $1\leq i \leq N$ as
$$t_N(P_i)=a_{s,i}t_l(P_N)=a_{s,i}b_{l,N}.$$
\\
We see now how to invert $B$, supposing 
to know the matrices $B'$ and $\begin{pmatrix}I'\\ \cline{1-1} C'\\ \end{pmatrix}$ associated to $\mathbf{X}'$.
We first border $B'$ getting a new matrix $B''$, by means of the remark on $\sigma$-antecents and $\sigma$-values stated above.
In particular, if $B'=(b'_{ij})$,  $1 \leq i,j \leq N-1$ then $B''=(b''_{h,m})$, $1 \leq h,m \leq N$ with 
$$
b''_{hm}=
\begin{cases}
 b'_{hm} \textrm{ if } 1 \leq h,m \leq N-1\\
 1 \textrm{ if } m=N, \, h=1\\
 a_{s(h),m}b''_{l(h)m} \textrm{ if } m=N,\, 2\leq h \leq N-1,\\
 a_{s(h),N}b''_{l(h)m} \textrm{ if }  h=N,\, 1\leq m \leq N 
\end{cases}
$$
where $s(h),l(h) $ are such that $t_h=x_{s(h)}t_{l(h)}$ so $P_{l(h)}$ is the last $\sigma$-antecedent of $P_h$
 and $s(h)$ the last $\sigma$-value. Then we
border $C'=(c'_{ij})$,  $1 \leq i,j \leq N-1$ with the last row and column of the  $N\times N$ identity matrix and we get 
$C''=(c''_{h,m})$, $1 \leq h,m \leq N$, with 
$$
c''_{hm}=
\begin{cases}
 c'_{hm} \textrm{ if } 1 \leq h,m \leq N-1\\
 0 \textrm{ if } h=N,\, 1\leq m \leq N-1 \textrm{ or } m=N,\, 1\leq h \leq N-1,\\
 1 \textrm{ if } h=m=N. 
\end{cases}
$$
In order to border $I'=(i'_{lj})$, $1\leq l,j\leq N-1$,  and get $I''=(i''_{hk})$, $1\leq h,k\leq N$,
we perform the following multiplication, which is a reformulation of (\ref{fBe})
\begin{equation}\label{iBC}
 (i''_{N1}....i''_{N,N-1})=(b''_{N1}....b''_{N,N-1})C'
\end{equation}
getting the terms $1,...,N-1$ of the last row of $I''$ (to be added on the bottom of $I'$); the last column of $I''$
is exactly the $N$-th column of   $B''$, in other words
$$
i''_{hk}=
\begin{cases}
 i'_{hk} \textrm{ if } 1 \leq h,k \leq N-1\\
 \sum_{j=1}^{N-1}b_{Nj}c'_{jk}  \textrm{ if } h=N,\, 1\leq k \leq N-1\\
b''_{hk} \textrm{ if } k=N,\, 1 \leq h \leq N.  
\end{cases}
$$
This way, we have  $\begin{pmatrix}I''\\ \cline{1-1} C''\\ \end{pmatrix}$; we reduce the last column 
w.r.t. the previous ones, getting 
$$
i''_{hk}=
\begin{cases}
 i'_{hk} \textrm{ if } 1 \leq h,k \leq N-1\\
 \sum_{j=1}^{N-1} b_{Nj}c'_{jk}  \textrm{ if } h=N,\, 1\leq k \leq N-1\\
 0 \textrm{ if } k=N,\, 1 \leq h \leq N-1\\  
b''_{NN}-\sum_{j=1}^{N-1}i''_{jN}i''_{Nj}=b''_{NN}-\sum_{j=1}^{N-1}b''_{jN}i''_{Nj}\textrm{ if } h=k=N
 \end{cases}
$$
Then we use again column reduction in order to reduce the last row 
of $I''$, so that $I''$ is transformed to the $N\times N$ identity matrix $I$. 
More precisely we first divide the last column of   $\begin{pmatrix}I''\\ \cline{1-1} C''\\ \end{pmatrix}$ by $i''_{NN}$
so each $i''_{hN}$ is transformed into $i''_{hN}/i''_{NN}$ and each $c''_{hN}$ is transformed into $c''_{hN}/i''_{NN}$, for $1 \leq h \leq N$
(in particular, now $i''_{NN}$ becomes $1$). Then, for $1 \leq m \leq N-1$, we substract the $N$-th column of   $\begin{pmatrix}I''\\ \cline{1-1} C''\\ \end{pmatrix}$
multiplied by $i''_{Nm}$ from the $m$-th column of  $\begin{pmatrix}I''\\ \cline{1-1} C''\\ \end{pmatrix}$ and we substitute 
the result to the $m$-th column of  $\begin{pmatrix}I''\\ \cline{1-1} C''\\ \end{pmatrix}$, obtaining  
a new matrix $\begin{pmatrix}I\\ \cline{1-1} C\\ \end{pmatrix}$, where $I$ is the $N \times N$ identity matrix and 
$C=(c_{hm})$, $1 \leq h,m \leq N$ s.t.

$$
c_{hm}=
\begin{cases}
 c''_{hm}-i''_{Nm}c''_{hN} \textrm{ if } 1 \leq h \leq N,\, 1 \leq m \leq N-1\\
 c''_{hm} \textrm{ if } m=N,\, 1\leq h \leq N.
\end{cases}
$$
Setting $B:=B''$, we have $C=B^{-1}$.
\\
Now, we compute $D^{(h)}$ edging the analogous matrix $D'^{(h)}$ for $\mathbf{X}'$, 
using an analogous argument as for $B$. In particular  
 if $D'^{(h)}=(d'^{(h)}_{lj})$,  $1 \leq l,j \leq N-1$ then $D''^{(h)}=(d''^{(h)}_{k,m})$, $1 \leq k,m \leq N$ with 
$$
d''^{(h)}_{km}=
\begin{cases}
 d'^{(h)}_{km} \textrm{ if } 1 \leq k,m \leq N-1\\
 a_{s(N),m}d''^{(h)}_{l(N)m} \textrm{ if } k=N,\, 1\leq m \leq N\\
  a_{s(k),m}d''^{(h)}_{l(k)m}  \textrm{ if } m=N,\, 1\leq k \leq N-1
\end{cases}
$$
where $s(k),l(k) $ are such that $t_k=x_{s(k)}t_{l(k)}$ so $P_{l(k)}$ is the last $\sigma$-antecedent of $P_k$
 and $s(k)$ the last $\sigma$-value.
\\
Finally, if $A'_{x_h}:=\left(a_{l'j'}^{'(h)}\right)_{l'j'}, 1\leq h\leq n$, $1\leq l',j'\leq N-1$ are the Auzinger-Stetter matrices for $\ck[x_1,...,x_n]/I(\mathbf{X}')$, we can find those for $A$, 
i.e. $A_{x_h}:=\left(a_{lj}^{(h)}\right)_{lj}, 1\leq h\leq n$, $1\leq l,j\leq N$, by 
computing 

$$
a_{l,j}=
\begin{cases}
 a'_{l,j}+d_{l,N}c_{N,j} \textrm{ if } l,j=1,...,N-1\\
 \sum_{k=1}^N d_{Nk}c_{kj}\textrm{ if } l=N, j=1,...,N\\
 \sum_{k=1}^N d_{lk}c_{kN} \textrm{ if } j=N,\, l=1,...,N.
 \end{cases}
$$
 
\begin{example}\label{EsAuStet}
  Consider again the set $\mathbf{X}=\{P_1=(1,0),P_2=(0,1),P_3=(0,2)\}$ of example \ref{Witnessmat} and \ref{separatoriconfrontati}.
\\ 
It is clear that for $P_1$, $B=C=1$ and $D^{(1)}=(1)=A_{x}$, $D^{(2)}=(0)=A_{y}$. 
\\ Adding $P_2$, we have $B''=\begin{pmatrix} 1\; 1\\1 \;0 \end{pmatrix}$ and
$\begin{pmatrix}I''\\ \cline{1-1} C''\\ \end{pmatrix}= \begin{pmatrix} 1\; 1\\ 1\; 0\\ \cline{1-1} 1\;0\\ 0\;1 \\ \end{pmatrix}$. We perform Gauss column reduction on $\begin{pmatrix}I''\\ \cline{1-1} C''\\ \end{pmatrix}$ by exchanging the two  columns and then
substituting the second column by the column obtained subtracting the first to the second one and we get $ \begin{pmatrix} 1\;\; 0\\ 0\;\; 1\\ \cline{1-1} 0\;\;1\\ 1-1 \\ \end{pmatrix}$, so
for this step  $B=\begin{pmatrix} 1\; 1\\1 \;0 \end{pmatrix}$ and  $C= \begin{pmatrix}   0\;\;1\\ 1-1 \\ \end{pmatrix}$ and $C=B^{-1}$.\\
Moreover, $D^{(1)}= \begin{pmatrix}   1\;0\\ 1\;0 \\ \end{pmatrix}$,
$A_x= \begin{pmatrix}   0\;1\\ 0\;1 \\ \end{pmatrix}$, $D^{(2)}= \begin{pmatrix}   0\;1\\ 0\;0 \\ \end{pmatrix}$,
$A_y= \begin{pmatrix}   1-1\\ 0\;\;0 \\ \end{pmatrix}$.
\\
Then we add also $P_3$ and we edge $B$ getting $B''=\begin{pmatrix} 1\; 1\;1\\1 \;0\;0\\0\;1\;2 \end{pmatrix}$  and  $C''= \begin{pmatrix}   0\;\;1\;\;0\\ 1-1\;0\\0\;\; 0\;\;1 \\ \end{pmatrix}$.
\\ Since $(0\; 1\; 2)\cdot C'' =(1\,-1\;2)$, then 
$I''=\begin{pmatrix}   1\;\;0\;\;1\\ 0\;\;1\;\;0\\1-1\;\-2 \\ \end{pmatrix}$.
We Gauss reduce $\begin{pmatrix}I''\\ \cline{1-1} C''\\ \end{pmatrix}$ via :

 $$ \begin{pmatrix}  1\;\;0\;\;1\\ 0\;\;1\;\;0\\1-1\;\-2 \\ \cline{1-1}    0\;\;1\;\;0\\ 1-1\;0\\0\;\; 0\;\;1 \\ \end{pmatrix} \rightarrow
   \begin{pmatrix}  1\;\;0\;\;0\\ 0\;\;1\;\;0\\1-1\;1 \\ \cline{1-1}    0\;\;1\;\;0\\ 1-1-1\\0\;\; 0\;\;1\\ \end{pmatrix}
\rightarrow
      \begin{pmatrix}  1\;\;0\;\;0\\ 0\;\;1\;\;0\\1\;\;0\;\;1 \\ \cline{1-1}    0\;\;1\;\;0\\ 1-2-1\\0\;\; 1\;\;1 \\ \end{pmatrix}  \rightarrow 
     \begin{pmatrix} 1\;\; 0 \;\; 0\\ 0\;\; 1\;\; 0\\ 0\;\;0\;\; 1\\ \cline{1-1} 0\;\;1\;\; 0\\ 2-2-1\\-1\;1\;\;1 \\ \end{pmatrix},
   $$
so $C=B^{-1}= \begin{pmatrix}  0\;\;1\;\; 0\\ 2-2-1\\-1\;1\;\;1 \\ \end{pmatrix}$. \\Finally,  $D^{(1)}= \begin{pmatrix}   1\;0\;0\\ 
1\;0\;0 \\0\;0\;0\\ \end{pmatrix}$,
$A_x= \begin{pmatrix}   0\;1\;0\\  0\;1\;0\\ 0\;0\;0 \\ \end{pmatrix}$, $D^{(2)}= \begin{pmatrix}   0\;1\;2\\ 0\;0\;0 \\0\;1\;4 \end{pmatrix}$,
$A_y= \begin{pmatrix}   0\;\;0\;\;1\\ 0\;\;0\;\;0 \\ 2-2\;\;3\\ \end{pmatrix}$.

\end{example}

\appendix
\section{Pseudocode of the algorithm}\label{PseudoCode}
In this section, we display the pseudocode of the algorithm just explained; 
in order to do so, we suppose known the following subroutines:
\begin{itemize}
 \item \texttt{children}, which takes as input a trie and one of its nodes,  and returns the ordered list of its children nodes.
 \item \texttt{parent}:  which takes as input  a trie and one of its nodes, and returns the father node.
  \item \texttt{edges}:  which takes as input  a trie and one of its nodes, and returns the edges   from this nodes to its children.
  \item \texttt{level}:  which takes as input  a trie and one of its nodes, and returns the level of the node, i.e. its distance from the root.
 \item \texttt{siblings}:  which takes as input  a trie and one of its nodes, and returns the ordered list of its siblings.
 \item \texttt{NewChild}:  which takes as input  a trie, one of its nodes $v$ and two labels  $\{i\}$ and $a$; it creates a new child for $v$, with label $\{i\}$ and edge from $v$ to the new child labelled with $a$
 \item \texttt{flatten}: which takes as input an ordered list of lists and returns an ordered list of its elements, 
 removing all the nested brackets\footnote{ We will use it to read the points 
corresponding to some bar in the Bar Code. Indeed, a Bar Code is stored 
 as a list of lists (see section \ref{core}).}.
 \item \texttt{min}: which takes as input an increasingly ordered list of positive integers $L$ and a list of positive integers $S$ and  
 finds the minimal element of $S$ contained in $L$. It returns $0$ if the intersection is empty.
\end{itemize}
\begin{example}\label{subroutines}
 Given the final point trie $T=\mathfrak{T}({\bf X})$ of example \ref{BC Trie}\\
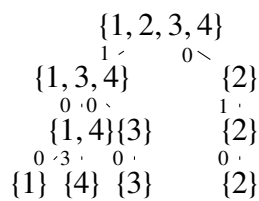
 \begin{minipage}{7cm}
 
 \begin{center}
 \begin{tikzpicture}[>=latex, scale=0.7] 
 \tikzstyle{cerchio}=[circle, draw=black,thick]
 \tikzstyle{cerchior}=[circle, draw=red,thick]
 \tikzstyle{rettangolo}=[rectangle, draw=black,thick]
\node (N_1) at (-0.5,5.75) {$\{1,2,3,4\}$};
\node (N_2) at (-2,4.75) {$\{1,3,4\}$};
\node (N_3) at (-2,3.75) {$\{1,4\}$};
\node (N_4) at (-3,2.75) {$\{1\}$};

\node (N_5) at (1,4.75) {$\{2\}$};
\node (N_6) at (1,3.75) {$\{2\}$};
\node (N_7) at (1,2.75) {$\{2\}$};

\node (N_8) at (-1,3.75) {$\{3\}$};
\node (N_9) at (-1,2.75) {$\{3\}$};
\node (N_10) at (-2,2.75) {$\{4\}$};

\draw [] (N_1) --node[left ]{$\scriptstyle{1}$} (N_2);
\draw [] (N_2) --node[left ]{$\scriptstyle{0}$} (N_3);
\draw [] (N_3) -- node[left ]{$\scriptstyle{0}$}(N_4);

\draw [] (N_1) --node[left ]{$\scriptstyle{0}$} (N_5) ;
\draw [] (N_5) --node[left ]{$\scriptstyle{1}$} (N_6);
\draw [] (N_6) -- node[left ]{$\scriptstyle{0}$}(N_7);

\draw [] (N_2) --node[left ]{$\scriptstyle{0}$} (N_8);
\draw [] (N_8) -- node[left ]{$\scriptstyle{0}$}(N_9);
\draw [] (N_3) -- node[left ]{$\scriptstyle{3}$}(N_10);

\end{tikzpicture}
\end{center}
 
 \end{minipage}
\hspace{0.5cm}
 \begin{minipage}{7cm}
  we have\\ \texttt{children}($T,\{1,3,4\}$)$=[\{1,4\},\{3\}]$, \\\texttt{parent}($T,\{1,4\}$)$=[\{1,3,4\}]$,\\
  \texttt{siblings}($T,\{1,4\}$)$=[\{3\}]$, \\  \texttt{edges}($T,\{1,4\}$)$=\{0,3\}$.\\
  \texttt{level}($T,\{1,4\}$)$=2$
 \end{minipage}\\
Applying NewChild($T,\{1,4\}$)($\{5\}, 2$), we get 
  
 
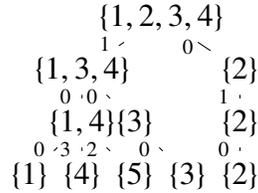
\begin{center}
 \begin{tikzpicture}[>=latex, scale=0.7] 
 \tikzstyle{cerchio}=[circle, draw=black,thick]
 \tikzstyle{cerchior}=[circle, draw=red,thick]
 \tikzstyle{rettangolo}=[rectangle, draw=black,thick]
\node (N_1) at (-0.5,5.75) {$\{1,2,3,4\}$};
\node (N_2) at (-2,4.75) {$\{1,3,4\}$};
\node (N_3) at (-2,3.75) {$\{1,4\}$};
\node (N_4) at (-3,2.75) {$\{1\}$};

\node (N_5) at (1,4.75) {$\{2\}$};
\node (N_6) at (1,3.75) {$\{2\}$};
\node (N_7) at (1,2.75) {$\{2\}$};

\node (N_8) at (-1,3.75) {$\{3\}$};
\node (N_11) at (-1,2.75) {$\{5\}$};

\node (N_9) at (0,2.75) {$\{3\}$};
\node (N_10) at (-2,2.75) {$\{4\}$};

\draw [] (N_1) --node[left ]{$\scriptstyle{1}$} (N_2);
\draw [] (N_2) --node[left ]{$\scriptstyle{0}$} (N_3);
\draw [] (N_3) -- node[left ]{$\scriptstyle{0}$}(N_4);
\draw [] (N_3) -- node[left ]{$\scriptstyle{2}$}(N_11);

\draw [] (N_1) --node[left ]{$\scriptstyle{0}$} (N_5) ;
\draw [] (N_5) --node[left ]{$\scriptstyle{1}$} (N_6);
\draw [] (N_6) -- node[left ]{$\scriptstyle{0}$}(N_7);

\draw [] (N_2) --node[left ]{$\scriptstyle{0}$} (N_8);
\draw [] (N_8) -- node[left ]{$\scriptstyle{0}$}(N_9);
\draw [] (N_3) -- node[left ]{$\scriptstyle{3}$}(N_10);

\end{tikzpicture}
\end{center}
Given a nested list $L=[[1,2,3],[[4,5],[6]]]$, we have \\ \texttt{flatten}($L$)$=[1,2,3,4,5,6]$,\\
\texttt{min}(\texttt{flatten}($L$), $\{7\}$)=$0$ and\\ \texttt{min}(\texttt{flatten}($L$), $\{2,3,4\}$)=$2$. 

\end{example}

We employ the poin trie construction as explained in 
section \ref{Lex}. More precisely, we exploit the procedure ExTrie, which takes as input a trie $T$ and a point $P_i$ and returns three inputs
\begin{itemize}
  \item the trie $T$, updated with the new point;
  \item the level $s$ in which the newly added path forks from the pre-existing trie; 
  \item the node $\mu$ at level $s$, s.t. its label is $V=\{i\}$; this is actually the node with minimum distance 
from the origin and s.t. its label is only $\{i\}$.
\end{itemize}
 \begin{algorithm}
\caption{ExTrie($T,P_i$): procedure extending a trie $T$ by adding a point $P_i$.}\label{ExTrie}
\begin{algorithmic}
\Require $T,P_i$
\Ensure $\rightarrow T,s,\mu$
\State root$(T)=$root$(T)\cup\{i\}$
\State $v=$root$(T)$
\State {\bf while } (level($T,v$)$<n$) {\bf do}
\State $f=0$ \Comment{Boolean variable, whose value, after the loop below, allows to decide whether to add a new child ($f=0$)
or not ($f=1$).}
\State $E=$edges($T,v$) \Comment{The edges from $v$ to its children}
\State {\bf for} $j=1$ to $\vert E \vert$
\If{$E[j]=P_i[\textrm{level}(T,v)+1]$}
\State $v=$children($v$)[$j$]$\cup \{i\}$
\State $f=1$ {\bf break}
\EndIf
\State {\bf end for}
\If{$f=0$}
\State $s=$level($T,v$)$+1$
\State $v=$NewChild($T,v$)($\{i\}, a_{s,i}$) \Comment{Create a new child for $v$ with label $\{i\}$ and edge from $v$ to the new child labelled with $a_{s,i}$}
\State $\mu=v$ \Comment{We keep track of the node in which the fork happens .}
\State {\bf for} $j=s+1$ to $n$ {\bf do}
\State $v=$NewChild($T,v$)($\{i\}, a_{j,i}$) \Comment{Create a new child for $v$ with label $\{i\}$ 
and edge from $v$ to the new child labelled with $a_{j,i}$.}
\State {\bf end for}
\EndIf
\State {\bf end while}
\State {\bf return } $T,s,\mu$
\end{algorithmic}
\end{algorithm}

The subroutine Fork, is used  to decide whether a path forks at some level $s$ of a trie $T$.
It takes as input a trie $T$, 
one of its nodes $v$, its level $s$ 
 and a set $S$ of identifiers for the points w.r.t. which we are looking for the fork.
 It returns $s$ if the fork happens and $0$ otherwise.
\begin{algorithm}
\caption{Fork($T,s, v,S$)}
\label{Fork}
\begin{algorithmic}[1]
\Require $T,s, v,S$
\Ensure $f$
\State $f=0$ \Comment{Variable whose value allows to decide whether the fork happens at level $s$ ($f = s$) or not ($f = 0$).}
\State $w=$ \texttt{parent}($T,v$)
\If{$w\cap S\neq \emptyset$}
\State $f=s$
\EndIf
\State {\bf return} $f$
\end{algorithmic}
\end{algorithm}

The subroutine $\sigma$ant is devoted to the quest for the $\sigma$-antecedent. It takes as input:
\begin{itemize}
\item  a trie $T$ and one of its nodes $v$;
\item a set $S$ of indices for the points among which we are looking for the $\sigma$-antecedent;
\end{itemize}
and it returns $l$, the index of the $\sigma$-antecedent.\\
\begin{algorithm}
\caption{$\sigma$-antecedent}\label{sant}
\begin{algorithmic}[1]
\Procedure{$\sigma$ant}{ $T,v,S$} $\rightarrow l$
\State $L=$\texttt{siblings}($T, v$);
\State $l=0$;
\State {\bf for } $j=\vert L\vert$ to $1$ {\bf do}
  \State  $l=$\texttt{min}($L[j]\cap S$);
   \If{$l \neq 0$}
    \State {\bf break};
    \EndIf
\State {\bf end for}
\State {\bf return} $l$
\EndProcedure
\end{algorithmic}
\end{algorithm}
The procedure NextB takes as input a Bar Code B, the $\sigma$-value $s$, the $\sigma$-antecedent's index $l$ and the matrix \footnote{Given a matrix $M=(m_{i,j})_{i\in I,j\in J}$ 
we denote its 
entries by $M[i][j]=m_{i,j}$.} $M$ and tests if we are in case a. or case b. (see section \ref{core}) for the bar\footnote{A given Bar Code $\cB$  is actually stored as a list of lists, so, a 
$i$-bar ($1 \leq i \leq n$) is identified by a list of indices (the e-list truncated to $s$). For example, if we take the Bar Code $\cB$ of example \ref{storeBC}, then $\cB^{(3)}_1=\cB[0]$ and 
$\cB^{2}_2=\cB[0,1]$.}
$B=\cB[M[l][1],...,M[l][s]]$. More precisely it tries to find $B'=\cB[M[l][1],...,M[l][s]+1]$; if we are in case a. it returns error, otherwise it returns $B'$.

\begin{algorithm}
\caption{Next bar}\label{NextB}
\begin{algorithmic}[1]
\Procedure{NextB}{ $\cB,s,l,M$} $\rightarrow B'$
\State {\bf return} $\cB[M[l][1],...,M[l][n-s+1]+1]$
\Comment{If such bar does not exist, the procedure returns an error}
\EndProcedure
\end{algorithmic}
\end{algorithm}

\begin{example}\label{UpdNextB}
  Consider the Bar code $\cB$ of examples 
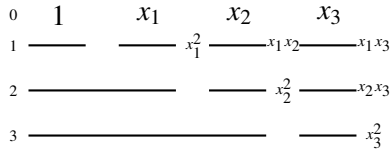
\ref{BCP} and \ref{elistExampleExp}.

 \begin{center}
\begin{tikzpicture}[scale=0.4]
\node at (-0.5,4) [] {${\scriptscriptstyle 0}$};
\node at (-0.5,0) [] {${\scriptscriptstyle 3}$};
\node at (-0.5,1.5) [] {${\scriptscriptstyle 2}$};
\node at (-0.5,3) [] {${\scriptscriptstyle 1}$};
 \draw [thick] (0,0) -- (7.9,0);
 \draw [thick] (9,0) -- (10.9,0);
 \node at (11.5,0) [] {${\scriptscriptstyle
x_3^2}$};
 \draw [thick] (0,1.5) -- (4.9,1.5);
 \draw [thick] (6,1.5) -- (7.9,1.5);
 \node at (8.5,1.5) [] {${\scriptscriptstyle
x_2^2}$};
 \draw [thick] (9,1.5) -- (10.9,1.5);
 \node at (11.5,1.5) [] {${\scriptscriptstyle
x_2x_3}$};
 \draw [thick] (0,3.0) -- (1.9,3.0);
 \draw [thick] (3,3.0) -- (4.9,3.0);
 \node at (5.5,3.0) [] {${\scriptscriptstyle
x_1^2}$};
 \draw [thick] (6,3.0) -- (7.9,3.0);
 \node at (8.5,3.0) [] {${\scriptscriptstyle
x_1x_2}$};
 \draw [thick] (9,3.0) -- (10.9,3.0);
 \node at (11.5,3.0) [] {${\scriptscriptstyle
x_1x_3}$};
 \node at (1,4.0) [] {\small $1$};
 \node at (4,4.0) [] {\small $x_1$};
 \node at (7,4.0) [] {\small $x_2$};
 \node at (10,4.0) [] {\small $x_3$};
\end{tikzpicture}
\end{center}
This is stored as \\
$\cB=\Biggl[\biggl[ \Bigl[[1],[3]\Bigr],\Bigl[[2]\Bigr]     \biggr],   \biggl[\Bigl[ [4]    \Bigr]\biggr] \Biggr].$\\
If we run NextB($\cB, 2, 3, M$), with
\[
M=\begin{bmatrix}
0& 0& 0\\
0& 1 & 0\\
0& 0 & 1
\end{bmatrix}
\]
it returns  $\cB[M[3][1],M[3][2]+1]=\cB[0,1]$, i.e. 
the second $2$-bar over the first $3$-bar, that is $\cB_2^{(2)}$ (the thick bar).

 
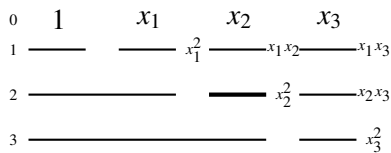
\begin{center}
\begin{tikzpicture}[scale=0.4]
\node at (-0.5,4) [] {${\scriptscriptstyle 0}$};
\node at (-0.5,0) [] {${\scriptscriptstyle 3}$};
\node at (-0.5,1.5) [] {${\scriptscriptstyle 2}$};
\node at (-0.5,3) [] {${\scriptscriptstyle 1}$};
 \draw [thick] (0,0) -- (7.9,0);
 \draw [thick] (9,0) -- (10.9,0);
 \node at (11.5,0) [] {${\scriptscriptstyle
x_3^2}$};
 \draw [thick] (0,1.5) -- (4.9,1.5);
 \draw [ultra thick] (6,1.5) -- (7.9,1.5);
 \node at (8.5,1.5) [] {${\scriptscriptstyle
x_2^2}$};
 \draw [thick] (9,1.5) -- (10.9,1.5);
 \node at (11.5,1.5) [] {${\scriptscriptstyle
x_2x_3}$};
 \draw [thick] (0,3.0) -- (1.9,3.0);
 \draw [thick] (3,3.0) -- (4.9,3.0);
 \node at (5.5,3.0) [] {${\scriptscriptstyle
x_1^2}$};
 \draw [thick] (6,3.0) -- (7.9,3.0);
 \node at (8.5,3.0) [] {${\scriptscriptstyle
x_1x_2}$};
 \draw [thick] (9,3.0) -- (10.9,3.0);
 \node at (11.5,3.0) [] {${\scriptscriptstyle
x_1x_3}$};
 \node at (1,4.0) [] {\small $1$};
 \node at (4,4.0) [] {\small $x_1$};
 \node at (7,4.0) [] {\small $x_2$};
 \node at (10,4.0) [] {\small $x_3$};
\end{tikzpicture}
\end{center}
If, instead, we run 
NextB($\cB, 2, 2, M$), we have  $\cB[M[2][1],M[2][2]+1]=\cB[0,2]$, so,  since there are no $2$-bars over $\cB^{(3)}_1$ and on the right of $\cB_2^{(2)}$, the procedure returns an error. 
\end{example}
The procedure Update, adds  new entries to the $N$-th row of a  matrix $M$ from the last column to that of the $\sigma$-value $s$. It takes as input the matrix $M$, the $\sigma$-value $s$, the 
$\sigma$-antecedent's identifier $l$ and the identifier of the new point, $N$, returning the updated $M$.

\begin{algorithm}
\caption{Update}\label{Update}
\begin{algorithmic}[1]
\Procedure{Update}{$M,s,l,N$} $\rightarrow M$
\State {\bf for } $i=1$ to $s-1$ {\bf do}
\State $M[N][n-i+1]=0$
\State {\bf end for}
\State  $M[N][n-s+1]=M[l][n-s+1]+1$
\State {\bf return} $M$
\EndProcedure
\end{algorithmic}
\end{algorithm}
We can display now the pseudocode of our algorithm. It takes as input a finite set of distinct points ${\bf X}$
and returns the matrix $M$, whose rows represent the terms corresponding to the elements of 
${\bf X}$, as explained in section \ref{core}.
\begin{algorithm}
\caption{Iterative Lex Game algorithm.}\label{IterativeLG}
\begin{algorithmic}[1]
\Procedure{IterLG}{$\underline{\mathbf{X}}$} $\rightarrow M$
\State $T=T(P_1)$ \Comment{Initialization}
\State $\cB=\cB(P_1)$
\State $M=[0,0,...,0]$
\State {\bf for } $i=2$ to $N$ {\bf do}
\State $T,f,v=ExTrie(T,P_i)$
\State {\bf for } $j=1$ to $n-f$ {\bf do}
\State $M[i][j]=0$
\State {\bf end for}
\State $S=\{P_1,...,P_{i-1}\}$
\State $s=f$
\State {\bf While} $s>0$ {\bf do}
\If{$f=s$}
\State  $l=\sigma ant(T, v, S)$
\State $B'=$NextB$(\cB,s,l,M)$
\If{$B'=$ error}
\State $\cB[M[l][1],...,M[l][s]+1,0..,0]=i$
\State $M =$Update$(M, s, l, i)$
\State $s=0$
\Else \State $M[i][n-s+1]=M[l][n-s+1]+1$
\State $S=$\texttt{flatten}$(B')$
\State $s=s-1$
\State $v=$\texttt{parent}$(v)$
\State $f=Fork(T,s, v, S)$
\EndIf
\Else
\State $M[i][n-s+1]=0$
\State $s=s-1$
\State $v=$\texttt{parent}$(v)$
\State $f=Fork(T,s, v, S)$
\EndIf
\State {\bf end while}
\State {\bf end for}
\State {\bf return} $M$
\EndProcedure
\end{algorithmic}
\end{algorithm}
\begin{Remark}\label{Msolounavolta}
  We point out that, for each row $1 \leq i\leq N$, each entry $M[i][j]$,
  $1 \leq j \leq n$ is written only once, with 
  \begin{itemize}
    \item $0$, in rows 4, 8, 27;
    \item $M[N][n-s+1]=M[l][n-s+1]+1$ in row  21;
    \item by the subroutine Update in row 18.
  \end{itemize}

\end{Remark}

\section{A (long) commented example }\label{AppA}
In this section we consider the set $\mathbf{X}=\{(1, 1, 2, 3), (1, 1, 2, 4), (1, 1, 2, 5), (1, 2, 1, 1), (1, 2, 1, 2), $ 
$ (1, 2, 2, 1), (1, 2, 2, 2), (3, 1, 1, 2), (3, 1, 2, 2), (3, 1, 2, 3), 
(3, 3, 1, 1), (3, 4, 1, 1),(3, 4, 1, 2)\}$, proposed for the first time by Gao-Rodrigues-Stroomer in \cite{GRS} and 
the ring $k[x_1,x_2,x_3,x_4]$, equipped with the lexicographical order induced by 
$1<x_1<x_2<x_3<x_4$.
\\
We start dealing with the point $P_1=(1,1,2,3)$, for which we know that $I(\{P_1\})=(x_1-1, x_2-1, x_3-2,x_4-3)$ and 
$\cN(I(\{P_1\}))=\{1\}$. We display the point trie and the Bar Code at this first step below:

\begin{minipage}{5cm}
 \begin{center}
\\
 that, represented as a list of lists, is $\left[\Biggl[\biggl[ \Bigl[[1] \Bigr]\biggr] \Biggr]\right].$
\end{center}
\end{minipage}\\
Then $t_1=1$ and
\[
M=%
\]

%%%%%%%%%%%%%%%%%%%%%%%%%%%%%%%%%%%%%%%%%
Adding $P_2=(1,1,2,4)$ to the point trie, we get

\begin{center}
 %
\end{center}
The fork happens in $s=4$ and the $\sigma$-antecedent is $P_l,$ for $l=1$, $B=\cB_1^{(4)}$ 

 \begin{center}
%
\end{center}

and there is no $4$-bar on the right of $B$, so we create $B'$
getting 
\begin{center}
%
\end{center}
that, represented as a list of lists, is $\left[\Biggl[\biggl[ \Bigl[[1] \Bigr]\biggr] \Biggr],\Biggl[\biggl[ \Bigl[[2] \Bigr]\biggr] \Biggr]\right].$

Then $t_2=x_4$ and
\[
M=%
\]

%%%%%%%%%%%%%%%%%%%%%%%%%%%%%%%%%%%%%%%%%
The insertion of $P_3=(1,1,2,5)$ in the point trie gives the following result 
\begin{center}
 %
\end{center}

The fork happens in $s=1$ and the $\sigma$-antecedent is $P_l,$ for $l=2$, so $B=\cB_2^{(4)}$ and again there is no $4$-bar on its right, so we create $B'$ and we get

\begin{center}
%
\end{center}
that, represented as a list of lists, is $\left[\Biggl[\biggl[ \Bigl[[1] \Bigr]\biggr] \Biggr],\Biggl[\biggl[ \Bigl[[2] \Bigr]\biggr] \Biggr],\Biggl[\biggl[ \Bigl[[3] \Bigr]\biggr] \Biggr]\right].$

Then $t_3=x_4^2$ and
\[
M=%
\]

%%%%%%%%%%%%%%%%%%%%%%%%%%%%%%%%%%%%%%%%%

We insert now $P_4=(1,2,1,1)$ in the point trie, getting 
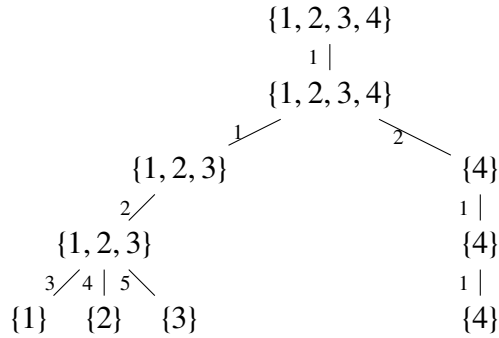
\begin{center}
 %
\end{center}
The fork happens in $s=2$ and 
$l=1$ so $B=\cB_1^{(2)}$

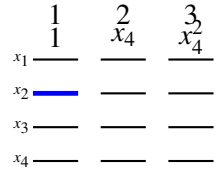
\begin{center}
%
\end{center}

and there is no $2$-bar on the right of $B$ lying over $\cB_1^{(4)},\cB_1^{(3)}$, so we create it and
we get 
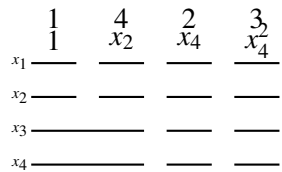
\begin{center}
%
\end{center}
that, represented as a list of lists, is $\left[\Biggl[\biggl[ \Bigl[[1] \Bigr], \Bigl[[4] \Bigr]\biggr] \Biggr],\Biggl[\biggl[ \Bigl[[2] \Bigr]\biggr] \Biggr],\Biggl[\biggl[ \Bigl[[3] \Bigr]\biggr] \Biggr]\right].$

Then $t_4=x_2$ and
\[
M=%
\]
%%%%%%%%%%%%%%%%%%%%%%%%%%%%%%%%%%%%%%%%%
The insertion of $P_5=(1,2,1,2)$ produces the following point trie:
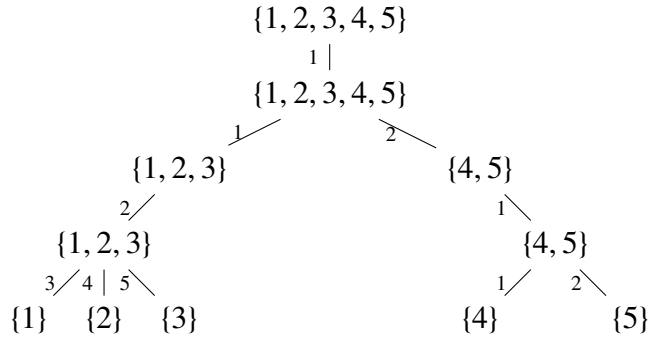
\begin{center}
 %
\end{center}
the fork happens in $s=4$ and the $\sigma$-antecedent is $P_l,$ for $l=4$, so $B=\cB_1^{(4)}$ and $B'=\cB_2^{(4)}$, thus $S=\{P_2\}$

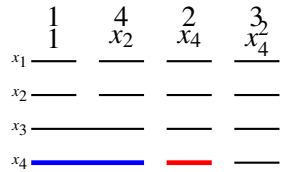
\begin{center}
%
\end{center}

\noindent Now, the fork with $S$ happens for $s=2$ and the $\sigma$-antecedent is $P_l,$ for $l=2$, so $B=\cB_3^{(2)}$ and, over $\cB_2^{(4)},\cB_2^{(3)}$ there is no 
$2$-bar following $B$
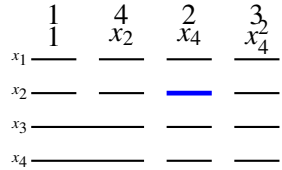
\begin{center}
%
\end{center}
\noindent so we create $B'$ getting 
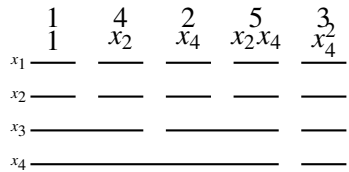
\begin{center}
%
\end{center}
\noindent that, represented as a list of lists, is $\left[\Biggl[\biggl[ \Bigl[[1] \Bigr], \Bigl[[4] \Bigr]\biggr] \Biggr],\Biggl[\biggl[ \Bigl[[2] \Bigr],\Bigl[[5] \Bigr]\biggr] \Biggr],\Biggl[\biggl[ \Bigl[[3] \Bigr]\biggr] \Biggr]\right].$

\noindent Then $t_5=x_2x_4$ and
\[
M=%
\]

%%%%%%%%%%%%%%%%%%%%%%%%%%%%%%%%%%%%%%%%%

\noindent Inserting $P_6=(1,2,2,1)$, we get the following point trie
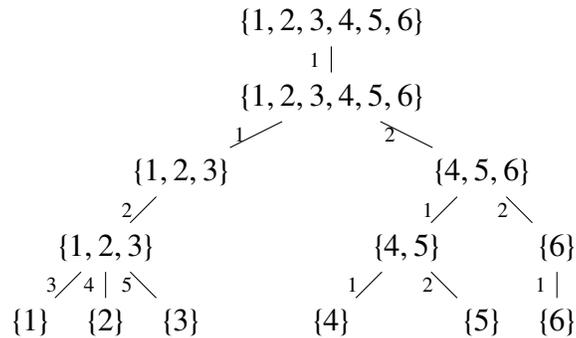
\begin{center}
 %
\end{center}
\noindent and, observing it, we see that the fork happens in $s=3$ and the $\sigma$-antecedent is $P_l,$ for $l=4$; $B=\cB_1^{(3)}$ and $B'$ has to be created
since there is no $3$-bar over $\cB_1^{(4)}$, on the right of $B$.

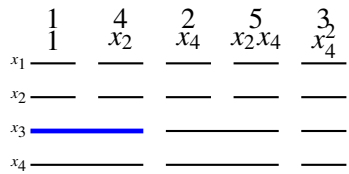
\begin{center}
%
\end{center}

\noindent We create then $B'$ getting

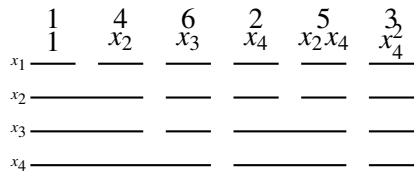
\begin{center}
%
\end{center}
\noindent that, represented as a list of lists, is $\left[\Biggl[\biggl[ \Bigl[[1] \Bigr], \Bigl[[4] \Bigr]\biggr], \biggl[ \Bigl[[6] \Bigr]\biggr] \Biggr],\Biggl[\biggl[ \Bigl[[2] \Bigr],\Bigl[[5] \Bigr]\biggr] \Biggr],\Biggl[\biggl[ \Bigl[[3] \Bigr]\biggr] \Biggr]\right].$

Then $t_6=x_3$ and
\[
M=%
\]

%%%%%%%%%%%%%%%%%%%%%%%%%%%%%%%%%%%%%%%%%

\noindent After inserting the point $P_7=(1,2,2,2)$, we have 
\begin{center}
 %
\end{center}
\noindent The fork happens for $s=4$ and the $\sigma$-antecedent is $P_l,$ for $l=6$ so $B=\cB_1^{(4)}$ and $B'=\cB_2^{(4)}$, so $S=\{P_2,P_5\}$

\begin{center}
%
\end{center}

\noindent The fork with $S$ happens for $s=3$ and the $\sigma$-antecedent is $P_l,$ for $l=5$ so $B=\cB_3^{(3)}$
\begin{center}
%
\end{center}

\noindent $B'$ is created since there is no $3$-bar on the right of $B$ and lying over $\cB_2^{(4)}$, so the new Bar Code is

\begin{center}
%
\end{center}
\noindent that, represented as a list of lists, is $\left[\Biggl[\biggl[ \Bigl[[1] \Bigr], \Bigl[[4] \Bigr]\biggr], \biggl[ \Bigl[[6] \Bigr]\biggr] \Biggr],\Biggl[\biggl[ \Bigl[[2] \Bigr],\Bigl[[5] \Bigr]\biggr]  , \biggl[ \Bigl[[7] \Bigr]\biggr]\Biggr],\Biggl[\biggl[ \Bigl[[3] \Bigr]\biggr]  \Biggr]\right].$

\noindent Then $t_7=x_3x_4$ and
\[
M=%
\]

%%%%%%%%%%%%%%%%%%%%%%%%%%%%%%%%%%%%%%%%%
\noindent Inserting in the point trie $P_8=(3,1,1,2)$, we get 
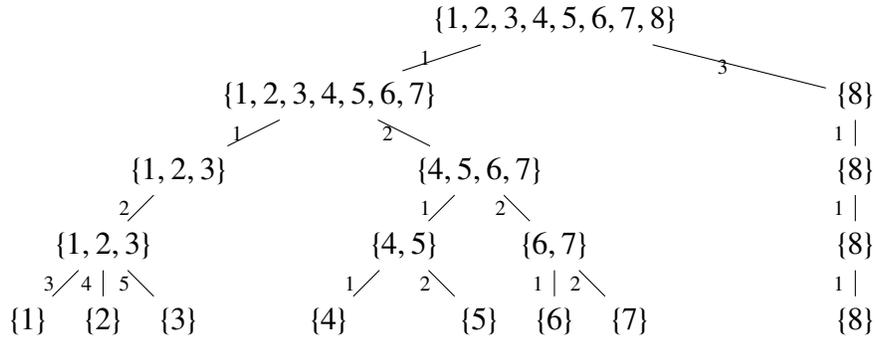
\begin{center}
 %
\end{center}
\noindent The fork happens for $s= 1$ and the $\sigma$-antecedent is $P_l,$ for $l=1$ so $B=\cB_1^{(1)}$ 

\begin{center}
%
\end{center}

\noindent and $B'$ has to be created, since no $1$-bar on the right of $B$ lies over
$\cB_1^{(4)},\cB_1^{(3)},\cB_1^{(2)}$, so we create it, getting 

\begin{center}
%
\end{center}
\noindent that, represented as a list of lists, is $\left[\Biggl[\biggl[ \Bigl[[1],[8] \Bigr], \Bigl[[4] \Bigr]\biggr], \biggl[ \Bigl[[6] \Bigr]\biggr] \Biggr],\Biggl[\biggl[ \Bigl[[2] \Bigr],\Bigl[[5] \Bigr]\biggr]  , \biggl[ \Bigl[[7] \Bigr]\biggr]\Biggr],\Biggl[\biggl[ \Bigl[[3] \Bigr]\biggr]  \Biggr]\right].$

Then $t_8=x_1$ and
\[
M=%
\]
%%%%%%%%%%%%%%%%%%%%%%%%%%%%%%%%%%%%%%%%%
\noindent Now we insert $P_9=(3,1,2,2)$, obtaining

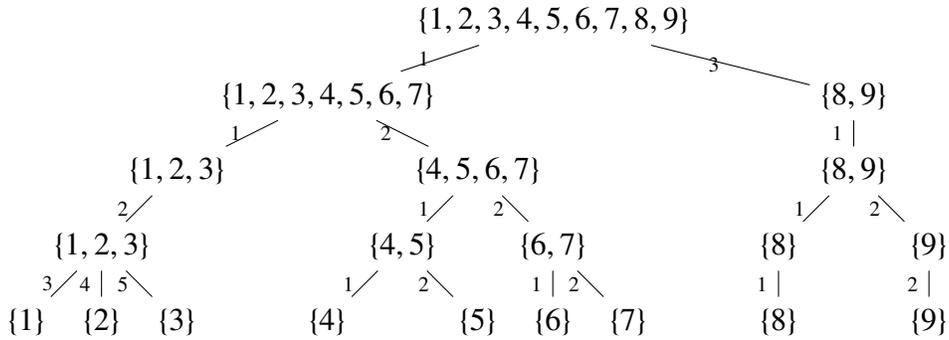
\begin{center}
 %
\end{center}

\noindent The fork is in $s=3$ and the $\sigma$-antecedent is $P_l,$ for $l=8$, so $B=\cB_1^{(3)}$, $B'=\cB_2^{(3)}$ and $S=\{P_6\}$

\begin{center}
%
\end{center}

\noindent Now the fork w.r.t. $S$ happens for $s=1$ and the $\sigma$-antecedent is $P_l,$ for $l=6$; $B=\cB_4^{(1)}$; we create $B'$ since there is no $1$-bar on the right 
of $B$, that lies over $\cB_1^{(4)},\cB_2^{(3)},\cB_3^{(2)}$ and we get

\begin{center}
%
\end{center}
\noindent that, represented as a list of lists, is $\left[\Biggl[\biggl[ \Bigl[[1],[8] \Bigr], \Bigl[[4] \Bigr]\biggr], \biggl[ \Bigl[[6],[9] \Bigr]\biggr] \Biggr],\Biggl[\biggl[ \Bigl[[2] \Bigr],\Bigl[[5] \Bigr]\biggr]  , \biggl[ \Bigl[[7] \Bigr]\biggr]\Biggr],\Biggl[\biggl[ \Bigl[[3] \Bigr]\biggr]  \Biggr]\right].$

Then $t_9=x_1x_3$ and
\[
M=%
\]

%\%%%%%%%%%%%%%%%%%%%%%%%%%%%%%%%%%%%%%%%%
\noindent We insert $P_{10}=(3,1,2,3)$, in the trie and we get 
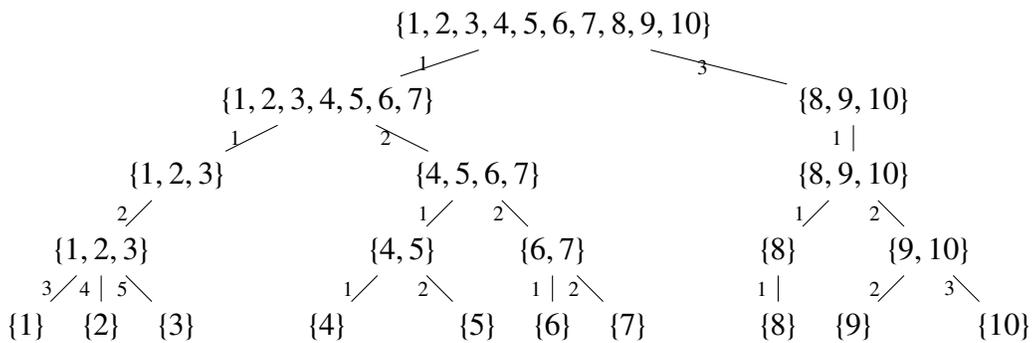
\begin{center}
 %
\end{center}

\noindent The fork happens for $s=4$ and the $\sigma$-antecedent is $P_l,$ for $l=9$, so $B=\cB_1^{(4)}$ and $B'=\cB_2^{(4)}$

\begin{center}
%
\end{center}

\noindent Then, $S=\{P_2,P_5,P_7\}$ and the fork with them is in $s=1$ and the $\sigma$-antecedent is $P_l,$ for $l=2$, so $B=\cB_6^{(1)}$ and $B'$
has to be created. We get 
\begin{center}
%
\end{center}
\noindent that, represented as a list of lists, is 
$$\left[\Biggl[\biggl[ \Bigl[[1],[8] \Bigr], \Bigl[[4] \Bigr]\biggr], \biggl[ \Bigl[[6],[9] \Bigr]\biggr] \Biggr],\Biggl[\biggl[ \Bigl[[2],[10] \Bigr],\Bigl[[5] \Bigr]\biggr]  , \biggl[ \Bigl[[7] 
\Bigr]\biggr]\Biggr],\Biggl[\biggl[ \Bigl[[3] \Bigr]\biggr]  \Biggr]\right].$$
\noindent Then $t_{10}=x_1x_4$ and
\[
M=%
\]

%%%%%%%%%%%%%%%%%%%%%%%%%%%%%%%%%%%%%%%%%
\noindent We insert now $P_{11}=(3,3,1,1)$ in the point trie
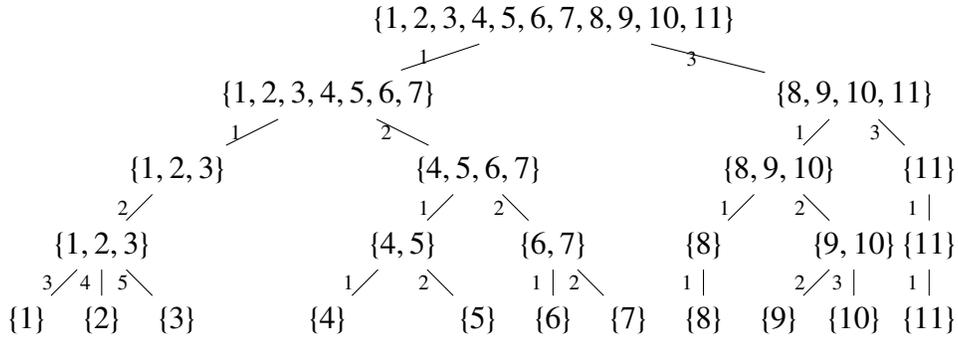
\begin{center}
 %
\end{center}
\noindent The fork happens in $s=2$ and the $\sigma$-antecedent is $P_l,$ for $l=8$, so $B=\cB_1^{(2)}$ and $B'=\cB_2^{(2)}$

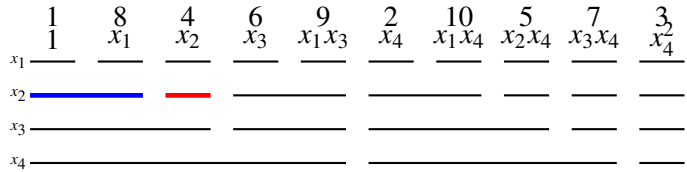
\begin{center}
%
\end{center}
\noindent Now, $S=\{P_4\}$; the fork happens in $s=1$ and the $\sigma$-antecedent is $P_l,$ for $l=4$, so $B=\cB_3^{(1)}$ and $B'$ 
has to be created, getting
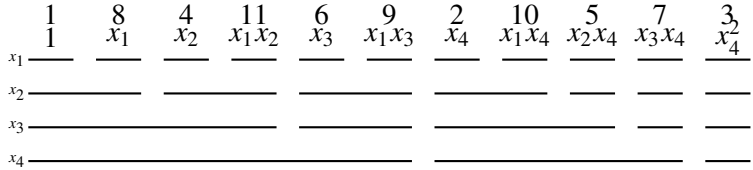
\begin{center}
%
\end{center}
\noindent that, represented as a list of lists, is 
$$\left[\Biggl[\biggl[ \Bigl[[1],[8] \Bigr], \Bigl[[4],[11] \Bigr]\biggr], \biggl[ \Bigl[[6],[9] \Bigr]\biggr] \Biggr],\Biggl[\biggl[ \Bigl[[2],[10] \Bigr],\Bigl[[5] \Bigr]\biggr]  , \biggl[ 
\Bigl[[7] \Bigr]\biggr]\Biggr],\Biggl[\biggl[ \Bigl[[3] \Bigr]\biggr]  \Biggr]\right].$$

\noindent Then $t_{11}=x_1x_2$ and
\[
M=%
\]

%%%%%%%%%%%%%%%%%%%%%%%%%%%%%%%%%%%%%%
\noindent Inserting $P_{12}=(3,4,1,1)$, we get the point trie that follows
\begin{center}
 %
\end{center}
\noindent The fork happens for $s=2$ and the $\sigma$-antecedent is $P_l,$ for $l=11$, indeed the 
rightmost sibling of the node $v_{2,5}$ s.t. $V_{2,5}=\{12\}$ is actually
$v_{2,4}$ and $V_{2,4}=\{11\}$. Then, $11$ is the leftmost element in 
 $V_{2,4}=\{11\}$.\\
We have $B=\cB_2^{(2)}$ and $B'$ has to be created, getting
\begin{center}
%
\end{center}
\noindent that, represented as a list of lists, is $$\left[\Biggl[\biggl[ \Bigl[[1],[8] \Bigr], \Bigl[[4],[11] \Bigr], \Bigl[[12] \Bigr]
\biggr], \biggl[ \Bigl[[6],[9] \Bigr]\biggr] \Biggr],\Biggl[\biggl[ \Bigl[[2],[10] \Bigr],\Bigl[[5] \Bigr]\biggr]  , \biggl[ \Bigl[[7] \Bigr]\biggr]\Biggr],\Biggl[\biggl[ \Bigl[[3] \Bigr]\biggr]  \Biggr]\right].$$

\noindent Then $t_{12}=x_2^2$ and
\[
M=%
\]

%%%%%%%%%%%%%%%%%%%%%%%%%%%
\noindent We finally insert $P_{13}=(3,4,1,2)$ and the final point trie is
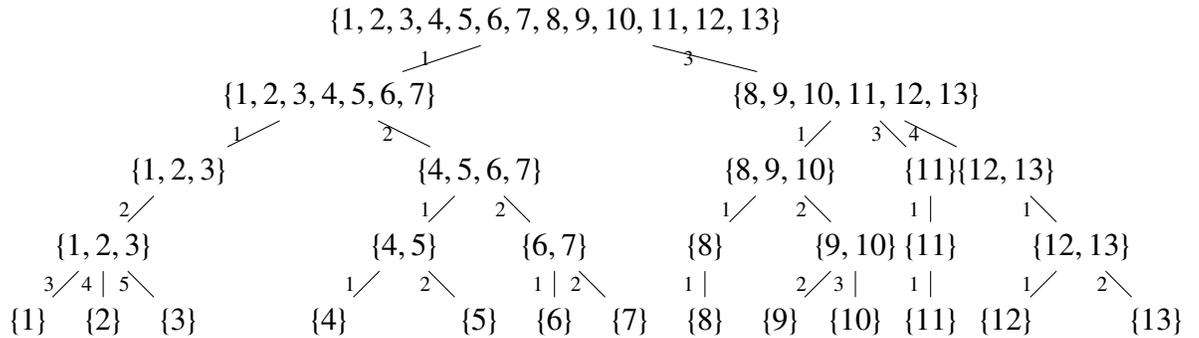
\begin{center}
 %
\end{center}
\noindent The fork happens for $s=4$ and the $\sigma$-antecedent is $P_l,$ for $l=12$, so $B=\cB_1^{(4)}$ 
 and $B'=\cB_2^{(4)}$ and $S=\{P_2,P_{10},P_5,P_7\}$

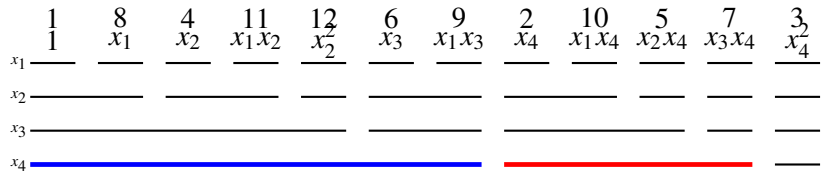
\begin{center}
%
\end{center}
\noindent The fork happens for $s=2$ and the $\sigma$-antecedent is $P_l,$ for $l=10$.
This happens since $S=\{P_2,P_{10},P_5,P_7\}$, so the rightmost sibling s.t. 
its label has nonempty intersection with $S$ is actually 
$v_{2,3}$, with $V_{2,3}=\{8,9,10\}$ and $10$ is the leftmost element in the intersection.
\\
We have $B=\cB_5^{(2)}$, $B'=\cB_6^{(2)}$, $S=\{P_5\}$

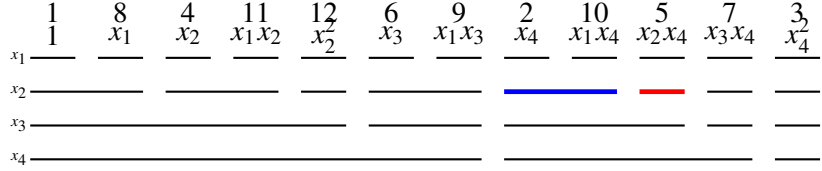
\begin{center}
%
\end{center}

\noindent The last fork is in $s=1$ and the $\sigma$-antecedent is $P_l,$ for $l=5$; $B=\cB_{10}^{(1)}$ and $B'$ is created.
\\
The final Bar code is then
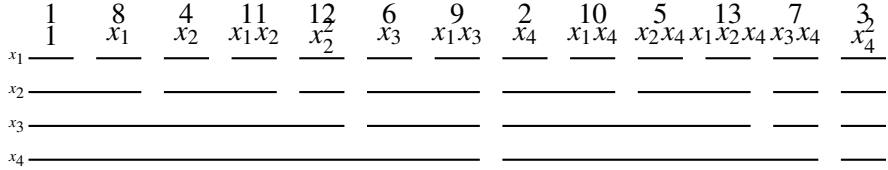
\begin{center}
%
\end{center}
\noindent that, represented as a list of lists, is $$\left[\Biggl[\biggl[ \Bigl[[1],[8] \Bigr], \Bigl[[4],[11] \Bigr], \Bigl[[12] \Bigr]
\biggr], \biggl[ \Bigl[[6],[9] \Bigr]\biggr] \Biggr],\Biggl[\biggl[ \Bigl[[2],[10] \Bigr],\Bigl[[5] ,[13]\Bigr]\biggr]  , \biggl[ \Bigl[[7] \Bigr]\biggr]\Biggr],\Biggl[\biggl[ \Bigl[[3] \Bigr]\biggr]  \Biggr]\right].$$

\noindent Then $t_{13}=x_1x_2x_4$ and
\[
M=%
\]

\end{document}